\documentclass[11pt, reqno]{amsart}

\usepackage[OT2, T1]{fontenc}

\usepackage[usenames,dvipsnames]{color}

\usepackage{url}
\usepackage{amsmath}
\usepackage{array}
\usepackage{graphicx}
\usepackage{amssymb}
\usepackage{amsthm}
\definecolor{link}{RGB}{11,0,128}
\usepackage[colorlinks=true,citecolor=NavyBlue,linkcolor=Bittersweet,urlcolor=link]{hyperref}
\usepackage{paralist}
\usepackage{colonequals}		% for nice :=
\usepackage{color}
\usepackage{mathabx}		% for \widec
\usepackage{amscd}
\usepackage[all,cmtip]{xy}	% for commutative diagrams
\usepackage{sseq}			% for spectral sequences
\usepackage{verbatim}
\usepackage{sverb}
\usepackage{fancyvrb}
\usepackage{parskip}
\usepackage{microtype}
\usepackage[in]{fullpage}
\usepackage{graphicx}
\usepackage{mathrsfs} 			% for \mathscr (script letters)
\usepackage{cleveref}
\usepackage{enumitem}
\usepackage{moreenum}			% for enumeration via Greek letters
\usepackage[alphabetic,lite]{amsrefs} 	% for bibliography
\usepackage{stmaryrd}	% this is for \llbracket and \rrbracket
\usepackage[foot]{amsaddr}
\usepackage{subfiles}
\usepackage{ragged2e}
\usepackage{relsize}

\DeclareSymbolFont{cyrletters}{OT2}{wncyr}{m}{n}
\DeclareMathSymbol{\Sha}{\mathalpha}{cyrletters}{"58}

\newcommand{\gL}{\lambda}

\newcommand{\bF}{\mathbb{F}}
\newcommand{\bG}{\mathbb{G}}

\newcommand{\bQ}{\mathbb{Q}}

\newcommand{\bT}{\mathbb{T}}

\newcommand{\bV}{\mathbb{V}}
\newcommand{\bZ}{\mathbb{Z}}

\newcommand{\cA}{\mathcal{A}}
\newcommand{\cB}{\mathcal{B}}
\newcommand{\cC}{\mathcal{C}}

\newcommand{\cE}{\mathcal{E}}

\newcommand{\cJ}{\mathcal{J}}
\newcommand{\cK}{\mathcal{K}}

\newcommand{\cO}{\mathcal{O}}

\newcommand{\cQ}{\mathcal{Q}}

\newcommand{\fm}{\mathfrak{m}}
\newcommand{\fn}{\mathfrak{n}}

\newcommand{\sE}{\mathscr{E}}

\newcommand{\sI}{\mathscr{I}}

\newcommand{\sO}{\mathscr{O}}

\newcommand{\ra}{\rightarrow}

\newcommand{\xra}{\xrightarrow}

\newcommand{\hra}{\hookrightarrow}

\newcommand{\wt}{\widetilde}
\newcommand{\wh}{\widehat}

\newcommand{\pr}{^{\prime}}

\newcommand{\ce}{\colonequals}
\newcommand{\ov}{\overline}

\newcommand{\sm}{\mathrm{sm}}

\renewcommand{\b}{\textbf}
\newcommand{\surjects}{\twoheadrightarrow}

\newcommand{\tensor}{\otimes} 		% binary tensor product
\newcommand{\isomto}{\overset{\sim}{\longrightarrow}}

\newcommand{\dR}{{\mathrm{dR}}}		% de Rham (mainly used in subscripts)
\newcommand{\et}{\mathrm{\acute{e}t}}	% for etale cohomology (mainly used in subscripts)
\newcommand{\llb}{\llbracket}		% [[
\newcommand{\rrb}{\rrbracket}		% ]]
\newcommand{\sh}{\mathrm{sh}}		% strict henselization (mainly used in superscripts)
\renewcommand{\i}{^{-1}}

\newcommand{\Div}{{\mathrm{Div}}}	% Cartier divisor associated to a generically exact complex

\providecommand{\In}[1]{\left\langle#1\right\rangle}
\providecommand{\p}[1]{\left(#1\right)}
\providecommand{\up}[1]{{\upshape(}#1{\upshape)}}
\providecommand{\uref}[1]{{\upshape\ref{#1}}}
\providecommand{\uS}{{\upshape\S}}
\providecommand{\ucolon}{{\upshape:} }
\providecommand{\uscolon}{{\upshape;} }

\providecommand{\f}[2]{\frac{#1}{#2}}
\DeclareMathOperator{\Ker}{Ker}			% Kernel
\DeclareMathOperator{\Coker}{Coker}		% Cokernel
\DeclareMathOperator{\im}{Im}			% Imaginary part
\DeclareMathOperator{\Spec}{Spec}		% Spectrum of a ring
\DeclareMathOperator{\Hom}{Hom}			% Set of arrows between two object
\DeclareMathOperator{\Char}{char}		% Characteristic of a field
\DeclareMathOperator{\Ext}{Ext}			% Derived functors of Hom
\DeclareMathOperator{\Extrig}{Extrig}			% Rigidified extensions
\DeclareMathOperator{\Gal}{Gal}	% Galois group
\DeclareMathOperator{\ord}{ord}	% order
\DeclareMathOperator{\GL}{GL}		% The general linear group
\DeclareMathOperator{\End}{End}		% The algebra of endomorphisms
\DeclareMathOperator{\Lie}{Lie}		% Lie algebra
\DeclareMathOperator{\length}{length}		% length of a module
\DeclareMathOperator{\Pic}{Pic}		% Picard group
\DeclareMathOperator{\Frob}{Frob}		% Frobenius
\DeclareMathOperator{\Fil}{Fil}			% filtration
\newcommand{\ba}{\begin{aligned}}
\newcommand{\ea}{\end{aligned}}
\newcommand{\be}{\begin{equation}}
\newcommand{\ee}{\end{equation}}
\newcommand{\pf}{\begin{proof}}
\newcommand{\bpf}{\begin{proof}}
\newcommand{\epf}{\end{proof}}
\newcommand{\bthm}{\begin{thm}}
\newcommand{\ethm}{\end{thm}}
\newcommand{\bthmt}{\begin{thm-tweak}}
\newcommand{\ethmt}{\end{thm-tweak}}
\newcommand{\bprop}{\begin{prop}}
\newcommand{\eprop}{\end{prop}}
\newcommand{\bcor}{\begin{cor}}
\newcommand{\ecor}{\end{cor}}
\newcommand{\bcort}{\begin{cor-tweak}}
\newcommand{\ecort}{\end{cor-tweak}}
\newcommand{\bconjt}{\begin{conj-tweak}}
\newcommand{\econjt}{\end{conj-tweak}}
\newcommand{\brem}{\begin{rem}}
\newcommand{\erem}{\end{rem}}

\newcommand{\brems}{\begin{rems} \hfill \begin{enumerate}[label=\b{\thesubsection.},ref=\thesubsection]}
\newcommand{\bremstweak}{\begin{rems-tweak} \hfill \begin{enumerate}[label=\b{\thesubsection.},ref=\thesubsection]}
\newcommand{\bremst}{\begin{rems-tweak} \hfill \begin{enumerate}[label=\b{\thesubsection.},ref=\thesubsection]}
\newcommand{\remi}{\addtocounter{subsection}{1} \item}

\newcommand{\remit}{\addtocounter{subsection}{1} \item}
\newcommand{\erems}{\end{enumerate} \end{rems}}
\newcommand{\eremstweak}{\end{enumerate} \end{rems-tweak}}
\newcommand{\eremst}{\end{enumerate} \end{rems-tweak}}
\newcommand{\blem}{\begin{lemma}}
\newcommand{\elem}{\end{lemma}}
\newcommand{\blemt}{\begin{lemma-tweak}}
\newcommand{\elemt}{\end{lemma-tweak}}
\newcommand{\bconj}{\begin{conj}}
\newcommand{\econj}{\end{conj}}
\newcommand{\bprob}{\begin{Problem}}
\newcommand{\eprob}{\end{Problem}}

\newcommand{\bq}{\begin{q}}
\newcommand{\eq}{\end{q}}
\newcommand{\benum}{\begin{enumerate}[label={{\upshape(\alph*)}}]}
\newcommand{\benuma}{\begin{enumerate}[label={{\upshape(\arabic*)}}]}
\newcommand{\benumr}{\begin{enumerate}[label={{\upshape(\roman*)}}]}
\newcommand{\eenum}{\end{enumerate}}
\newcommand{\bc}{\begin{comment}}
\newcommand{\ec}{\end{comment}}
\newcommand{\bd}{\begin{defn}}
\newcommand{\ed}{\end{defn}}
\newcommand{\forg}{\mathrm{forg}}
\newcommand{\quot}{\mathrm{quot}}

\newcommand{\beg}{\begin{eg}}
\newcommand{\eeg}{\end{eg}}

\newcommand{\bcl}{\begin{claim}}
\newcommand{\ecl}{\end{claim}}
\newcommand{\lab}{\label}

\newcommand{\q}{\quad}
\newcommand{\qq}{\quad\quad}
\newcommand{\qqq}{\quad\quad\quad}

\newcommand{\tst}{\textstyle}

\newcommand{\val}{\mathrm{val}}
\newcommand{\Inf}{\mathrm{inf}}
\newcommand*{\QED}{\hfill\ensuremath{\qed}}

\theoremstyle{plain}
\newtheorem{thm}[subsection]{Theorem}
\Crefname{thm}{Theorem}{Theorems}

\Crefname{rethm}{Theorem}{Theorem}
\newtheorem{prop}[subsection]{Proposition}
\Crefname{prop}{Proposition}{Propositions} 

\Crefname{Q}{Question}{Questions}
\newtheorem{Problem}[subsection]{Problem}
\Crefname{Problem}{Problem}{Problems}
\newtheorem{conj}[subsection]{Conjecture}
\Crefname{conj}{Conjecture}{Conjectures}
\newtheorem{cor}[subsection]{Corollary}
\Crefname{cor}{Corollary}{Corollaries}
\newtheorem{lemma}[subsection]{Lemma}

\Crefname{subprop}{Proposition}{Propositions}

\Crefname{subcor}{Corollary}{Corollaries}

\Crefname{sublem}{Lemma}{Lemmas}

\theoremstyle{remark}
\newtheorem{claim}[equation]{Claim}
\Crefname{claim}{Claim}{Claims}

\Crefname{subrem}{Remark}{Remarks}

\theoremstyle{definition}
\newtheorem{defn}[subsection]{Definition}
\Crefname{defn}{Definition}{Definitions}

\Crefname{conv}{Convention}{Conventions}
\newtheorem{eg}[subsection]{Example}
\Crefname{eg}{Example}{Examples}
\newtheorem{rem}[subsection]{Remark}
\Crefname{rem}{Remark}{Remarks}
\newtheorem*{rems}{Remarks}

\theoremstyle{plain}
\newtheorem{thm-tweak}[subsection]{Theorem}
\Crefname{thm-tweak}{Theorem}{Theorems}
\newtheorem{lemma-tweak}[subsection]{Lemma}
\Crefname{lemma-tweak}{Lemma}{Lemmas}
\newtheorem{cor-tweak}[subsection]{Corollary}
\Crefname{cor-tweak}{Corollary}{Corollaries}
\newtheorem{prop-tweak}[subsection]{Proposition}
\Crefname{prop-tweak}{Proposition}{Propositions} 
\newtheorem{conj-tweak}[subsection]{Conjecture}
\Crefname{conj-tweak}{Conjecture}{Conjectures} 

\theoremstyle{definition}
\newtheorem{defn-tweak}[subsection]{Definition}
\Crefname{defn-tweak}{Definition}{Definitions}
\newtheorem{eg-tweak}[subsection]{Example}
\Crefname{eg-tweak}{Example}{Examples}
\newtheorem*{rems-tweak}{Remarks}
\newtheorem{rem-tweak}[subsection]{Remark}
\Crefname{rem-tweak}{Remark}{Remarks}

\newtheoremstyle{subsection-tweak}
   {11pt}
   {3pt}%
   {}
   {}%
   {\bfseries}
   {}%
   {.5em}
   {\thmnumber{\@{#1}{}\@{#2}.}%
    \thmnote{~{\bfseries#3.}}}    
    
\theoremstyle{subsection-tweak}
\newtheorem{pp}[subsection]{}
\newcommand{\bpp}{\begin{pp}}
\newcommand{\epp}{\end{pp}}

\theoremstyle{subsection-tweak}
\newtheorem{pp-tweak}[subsection]{}

\numberwithin{equation}{subsection}

\makeatletter
\def\@tocline#1#2#3#4#5#6#7{
    \begingroup %\hyphenpenalty\@M
    \@ifempty{#4}{%
    }{%
    }%

    \parindent\z@ \leftskip#3\relax \advance\leftskip\@tempdima\relax
    #5\hskip-\@tempdima
      \ifcase #1
       \or\or \hskip 2em \or \hskip 1em \else \hskip 3em \fi%
      #6\nobreak\relax
    \dotfill\hbox to\@pnumwidth{\@tocpagenum{#7}}\par
    \nobreak
    \endgroup
  }
 \def\l@section{\@tocline{1}{0pt}{1pc}{}{}}

\renewcommand{\tocsection}[3]{%
  \indentlabel{\@ifnotempty{#2}{\makebox[1.3em][l]{%
    \ignorespaces#1 \bfseries{#2}.\hfill}}}\bfseries{#3}
    \vspace{1.5pt}}

\renewcommand{\tocsubsection}[3]{%
  \indentlabel{\@ifnotempty{#2}{\hspace*{-0.5em}\makebox[2.1em][l]{%
    \ignorespaces#1#2.\hfill}}}#3
    \vspace{1.5pt}}

\makeatother

\begin{document}

\title{The Manin constant in the semistable case}

\author{K\k{e}stutis \v{C}esnavi\v{c}ius}
\address{Laboratoire de Math\'{e}matiques d'Orsay, Univ.~Paris-Sud, CNRS, Universit\'{e} Paris-Saclay, 91405 Orsay, France}
\email{kestutis@math.u-psud.fr}
%\urladdr{http://math.berkeley.edu/~kestutis/}

\date{\today}

%------------------------------------------------------------------------------------------------

\begin{abstract} For an optimal modular parametrization $J_0(n) \surjects E$ of an elliptic curve $E$ over $\bQ$ of conductor $n$, Manin conjectured the agreement of two natural $\bZ$-lattices in the $\bQ$-vector space $H^0(E, \Omega^1)$. Multiple authors generalized his conjecture to higher dimensional newform quotients. We prove the Manin conjecture for semistable $E$, give counterexamples to all the proposed generalizations, and prove several semistable special cases of these generalizations. The proofs establish general relations between the integral $p$-adic \'{e}tale and de Rham cohomologies of abelian varieties over $p$-adic fields and exhibit a new exactness result for N\'{e}ron models. \end{abstract}

%-----------------------------------------------------------------------------------------------

\maketitle

\hypersetup{
    linktoc=page,     %set to all if you want both sections and subsections linked
}

\renewcommand*\contentsname{}
\q\\
\tableofcontents

%--------------------------------------------------------------------------End of metadata

\section{Introduction}

%\kestutis{Change $\subset$ to $\subseteq$ to emphasize that $H = \Gamma_0(n)$ is allowed?}

By the modularity theorem, every elliptic curve $E$ over $\bQ$ arises as a quotient 
\[
\pi \colon J_0(n) \surjects E
\]
of the modular Jacobian $J_0(n)$ with $n$ equal to be the conductor of $E$. In this situation there are two natural $\bZ$-lattices in the $\bQ$-vector space $\pi^*(H^0(E, \Omega^1))$ generated, respectively, by the pullback of a N\'{e}ron differential $\omega_E$ and by the $1$-form $f_E$ associated to the normalized newform determined by the isogeny class of $E$. The \emph{Manin constant} $c_\pi \in \bQ^\times$, defined by
\[
\pi^*(\omega_E) = c_\pi \cdot f_E,
\]
describes the difference between the two lattices and is the subject of the following conjecture of Manin (see \S\ref{conv} for a review of the terminology used in its formulation).

\bconjt[\cite{Man71}*{10.3}] \lab{manin-conj}
For an elliptic curve quotient 
\[
\pi\colon J_0(n) \surjects E
\]
that is \emph{new} in the sense that $n$ equals the conductor of $E$ and \emph{optimal} in the sense that $\Ker(\pi)$ is connected, the two $\bZ$-lattices in $\pi^*(H^0(E, \Omega^1))$ described above agree\uscolon in other words,
\[
c_\pi = \pm 1.
\]
\econjt

Substantial computational evidence supports the conjecture; for instance, Cremona proved that $c_\pi = \pm 1$ whenever $n \le 390000$, see \cite{Cre16} (and also \cite{ARS06}*{Thm.~2.6}). The main goal of this paper is to settle \Cref{manin-conj} for semistable $E$; more precisely, we prove the following result.

\bthmt[\Cref{manin-H-sst}] \lab{main-thm}
For an $n \in \bZ_{\ge 1}$, a subgroup $H \subset \GL_2(\wh{\bZ})$ with $\Gamma_1(n) \subset H \subset \Gamma_0(n)$, a new elliptic optimal quotient $\pi\colon J_H \surjects E$, and a prime $p$,
\[
\text{if} \q  p^2\nmid n, \q  \text{then}\q  \ord_p( c_\pi) = 0 \q \text{and $\pi$ induces a smooth morphism on N\'{e}ron models over $\bZ_p$.}
\]
 In particular, Conjecture \uref{manin-conj} holds in the case when $E$ is semistable \up{that is, when $n$ is squarefree}.
\ethmt

The proof of \Cref{main-thm} is relatively short, does not rely on any previously known cases, is uniform for all $p$, and is given in \S\ref{ell-curve}. The key idea is to translate multiplicity one results for differentials in characteristic $p$ into Hecke-freeness statements about the Lie algebra of the N\'{e}ron model of $J_0(n)$. This provides control on congruences between the ``$f_E$-isotypic'' and the ``(non-$f_E$)-isotypic'' parts of this Lie algebra, which combines with comparisons between the modular degree and the congruence number of $E$ to imply the $\Gamma_0(n)$ case of \Cref{main-thm}. The general case reduces to $\Gamma_0(n)$.

Previous results \cite{Maz78}*{Cor.~4.1}, \cite{AU96}*{Thm.~A and (ii) on p.~270}, \cite{ARS06}*{Thm.~2.7}, and \cite{Ces16g}*{Thm.~1.5} cover many special cases of \Cref{main-thm} (in particular, the case of an odd $p$). %---see \cite{Ces16g}*{\S1, esp.~Thm.~1.4 and Thm.~1.5} for a precise summary. 
Nevertheless, some semistable elliptic curves escape the net of these previous results: for instance, this happens for \href{http://www.lmfdb.org/EllipticCurve/Q/130/a/2}{$130.a2$}, \href{http://www.lmfdb.org/EllipticCurve/Q/4930/b/1}{$4930.b1$}, \href{http://www.lmfdb.org/EllipticCurve/Q/182410/a/1}{$182410.a1$}, and many others in \cite{LMFDB}. Beyond the semistable case, Edixhoven proved in \cite{Edi91}*{Thm.~3} that a new elliptic optimal quotient $\pi \colon J_0(n) \surjects E$ and a prime $p \ge 11$ for which $E_{\bQ_p}$ does not have potentially ordinary reduction of Kodaira type II, III, or IV satisfy $\ord_p(c_\pi) = 0$.

%Another consequence is derived in \S\ref{et-dR}: we show in \Cref{cl-im} that, in the setting of \Cref{main-thm}, if $H = \Gamma_0(n)$ and $p^2 \nmid n$, then $\pi^\vee$ induces a closed immersion between the N\'{e}ron models over $ \bZ_p$.

\bpp[Counterexamples to proposed generalizations] \lab{higher}
Generalizations of \Cref{manin-conj} to newform quotients of arbitrary dimension have been put forward by Conrad--Edixhoven--Stein \cite{CES03}*{Conj.~6.1.7}, Joyce \cite{Joy05}*{Conj.~2}, and Agashe--Ribet--Stein \cite{ARS06}*{Conj.~3.12}. These more general conjectures are supported by a handful of examples \cite{FLSSSW01}*{\S4.2}, \cite{ARS06}*{p.~624} and by the semistable case at odd $p$ that follows from exactness properties of N\'{e}ron models. %Part of the difficulty in extending the techniques of \cite{AU96} or \cite{Ces16f} to the higher dimensional case is that they rely on the exact structure of $\Ker(\pi \circ \pi^\vee)$ that is specific to elliptic curves.

We prove in \Cref{genl-counter} that all these generalizations of the Manin conjecture fail (at the prime $2$) for a $24$-dimensional optimal newform quotient of $J_0(431)$ and also for a $91$-dimensional optimal newform quotient of $J_0(2089)$.\footnote{Our counterexamples to the generalizations of the Manin conjecture rely on the correct functioning of the \cite{Sage} commands used in the proof of \Cref{genl-counter}.} On the positive side,  we prove a number of their semistable special cases in \Cref{genl-results}, which show that in the semistable setting the failure of the generalizations has to involve both the failure of mod $p$ multiplicity one for $J_0(n)$ and the failure of the maximality of the order $\cO_f$ determined by the newform $f$ (in the elliptic curve case $\cO_f$ is always $\bZ$).
\epp

\bpp[The interplay between integral \'{e}tale and de Rham cohomologies]
It is a natural idea that relations between the integral \'{e}tale and de Rham cohomologies of abelian varieties over $p$-adic fields are relevant for the Manin conjecture: one guesses that the role of the optimality assumption is to supply exactness on $H^1_\et(-, \wh{\bZ})$, whereas similar exactness on $H^1_\dR(-/\bZ)$ would be related to exactness on N\'{e}ron models implied by the Manin conjecture (more precisely, the Manin conjecture implies exactness on $\Fil^1(H^1_\dR(-/\bZ))$, that is, on $\Lie(-/\bZ)^*$). However, the required transfer of an integral \'{e}tale assumption to an integral de Rham conclusion has been problematic in the past, partly because comparisons of $p$-adic Hodge theory tend to fail integrally. The key novelty that underlies our approach is the idea that an arithmetic duality result of Raynaud, when reformulated as \Cref{Raynaud-re}, supplies an integral link between the two cohomology theories.

\Cref{Raynaud-re} is very robust for our purposes and is the backbone of \S\S\ref{et-dR}--\ref{contorni}. In addition to its role in the proofs of the results mentioned in \S\ref{higher}, it eventually supplies the exactness on $H^1_\dR(-/\bZ_p)$ in \Cref{cl-im} and \Cref{genl-results} (that is, in the cases of the Manin conjecture  and of its generalization proved in this paper) and it leads to cohomology specialization results in the spirit of \cite{BMS16} with no restrictions on the reduction type (see \Cref{BMS-eg}), to relations between torsion multiplicities in $J_0(n)$ and Gorenstein defects of Hecke algebras (see \Cref{Gor}), and to an exactness result for N\'{e}ron models equipped with a ``Hecke action'' (see \Cref{Ner-exact}).
\epp

\bpp[Notation and conventions] \lab{conv}
For an open subgroup $H \subset \GL_2(\wh{\bZ})$, we denote the level $H$ modular curve over $\bZ$ by $X_H$ (see \cite{Ces15a}*{6.1--6.3} for a review of $X_H$), and we denote the Jacobian of $(X_H)_\bQ$ by $J_H \ce \Pic^0_{(X_H)_\bQ/\bQ}$. For an $n \in \bZ_{\ge 1}$, we let $\Gamma_0(n)$ (resp.,~$\Gamma_1(n)$) denote the preimage of $\{ \p{\begin{smallmatrix} * & * \\ 0 & * \end{smallmatrix}} \} \subset \GL_2(\bZ/n\bZ)$ (resp.,~of $\{ \p{\begin{smallmatrix} 1 & * \\ 0 & * \end{smallmatrix}} \} \subset \GL_2(\bZ/n\bZ)$) in $\GL_2(\wh{\bZ})$, and we set $X_0(n) \ce X_{\Gamma_0(n)}$, etc. For an $H$ with $\Gamma_1(n) \subset H \subset \Gamma_0(n)$, a quotient $\pi \colon J_H \surjects A$ is \emph{optimal} (resp.,~\emph{new} or \emph{newform}) if $\Ker \pi$ is connected (resp.,~if up to isogeny $\pi$ arises from some newform $f_A$). We say that $f_A$ is \emph{normalized} if its $q$-expansion $a_1 q + a_2q^2 + \ldots$ at the cusp ``$\infty$'' has $a_1 = 1$. If $\pi$ is new and optimal, $f_A$ is normalized, and $\dim A = 1$, then we let $c_\pi$ be the \emph{Manin constant} defined by the equality $\pi^*(\omega_A) = c_\pi \cdot f_A$, where $\omega_A$ is a N\'{e}ron differential on $A$ and we have identified $f_A$ with its associated differential form on $J_H$. Up to a sign, $c_\pi$ does not depend on the choice of $\omega_A$.

For a commutative ring $\bT$, a maximal ideal $\fm \subset \bT$, and a $\bT$-module $M$, we let $M_\fm$ denote the $\fm$-adic completion of $M$ and we let $M[\fm^\infty]$ denote the submodule of the elements of $M$ killed by some power of $\fm$. We often consider $\bZ$-torsion free $\bT$ and $M$, for which we repeatedly abuse notation:
\be \lab{abuse}
M[e] \ce M \cap \Ker\p{e \colon M_\bQ \ra M_\bQ} \qq \text{for an idempotent} \qq e \in \bT_\bQ.
\ee
We let $\ord_p$ denote the $p$-adic valuation normalized by $\ord_p(p) = 1$. For a field $K$, we let $\ov{K}$ denote a fixed algebraic closure of $K$. For a commutative ring $R$ and a projective $R$-module $P$, we set $P^* \ce \Hom_R(P, R)$. For a smooth group scheme $G \ra S$, we let $G^0$ (resp.,~$\Lie G$) denote its relative identity component subfunctor (resp.,~its Lie algebra at the identity section), which in the situations below will always exist as a scheme. For a scheme $S$ and an $S$-scheme $X$, we let $X^\sm$ denote the smooth locus of $X$. When denoting structure sheaves or sheaves of K\"{a}hler differentials, we omit subscripts that may be inferred from the context. We let $(-)^\vee$ denote the dual of an abelian variety, or of a homomorphism of abelian varieties, or of a commutative finite locally free group scheme.

\epp

\subsection*{Acknowledgements}
I thank Kevin Buzzard, Frank Calegari, Brian Conrad, John Cremona, Naoki Imai, Bjorn Poonen, Michel Raynaud, Ken Ribet, Sug Woo Shin, and Preston Wake for helpful conversations or correspondence. I thank William Stein for sending me the Sage code used to compute the example in \cite{ARS12}*{Rem.~3.7}. I thank the referee for helpful comments and suggestions. I thank the Miller Institute for Basic Research in Science at the University of California Berkeley for support during the preparation of this article.

%%%%%%%%%%%%%%%%%%%%%%%%%%%%%%%%%%%%%%%

\section{The semistable case of the Manin conjecture} \lab{ell-curve}

The main goal of this section is to prove the semistable case of the Manin conjecture in \Cref{Manin-sst}. The path to this consists of the notational review in \S\ref{Hecke-def}, of parts \ref{H1f-i} and \ref{H1f-iii} of \Cref{H1-free}, and then of the discussion in between \S\ref{mod-ell-curve} and \Cref{cong-deg}. Modulo standard inputs from the literature, the overall argument is brief.

\bpp[The Hecke algebra $\bT$] \lab{Hecke-def}
Throughout \S\ref{ell-curve} we fix an $n \in \bZ_{\ge 1}$ and for primes $\ell \nmid n$ and $\ell \mid n$, respectively, we let $T_\ell$ and $U_\ell$ be the endomorphisms of $J_0(n)$ induced via ``Albanese functoriality'' by their namesake correspondences  (see \cite{MW84}*{Ch.~II, \S\S5.4--5.5}; in the notation there, we choose $T_{\ell^*}$ and $U_{\ell^*}$). The pullback action of $T_\ell$ and $U_\ell$ on $H^0(J_0(n), \Omega^1)$ agrees with their ``classical'' action on the space of weight $2$ cusp forms, see \cite{MW84}*{Ch.~II, \S5.8} (this is our reason for preferring the Albanese functoriality). We let 
\[
\bT \subset \End_\bQ(J_0(n))
\]
be the commutative $\bZ$-subalgebra generated by all the $T_\ell$ and $U_\ell$, so $\bT$ acts on various objects naturally attached to $J_0(n)$, e.g., 
\begin{itemize}
\item
on the N\'{e}ron model $\cJ$ over $\bZ$ of $J_0(n)$;

\item
on the tangent space $\Lie \cJ$ of $\cJ$ at the identity section and on the dual $H^0(\cJ, \Omega^1) = (\Lie \cJ)^*$.
\end{itemize}

If $p$ is a prime with $p^2 \nmid n$, then, by \cite{DR73}*{VI.6.7, VI.6.9}, the curve $X_0(n)_{\bZ_p}$ is semistable over $\bZ_p$, so that, by \cite{BLR90}*{9.7/2}, we have 
\[
(\Pic^0_{X_0(n)/\bZ})_{\bZ_p} \cong \cJ^0_{\bZ_p}.
\]
In particular, by \cite{BLR90}*{8.4/1 (a)}, \cite{Con00}*{Cor.~5.1.3}, and \cite{Con00}*{Thm.~B.4.1} (applied over $\bQ_p$), %(\cite{Con00}*{B.4.1} applied over $\bQ_p$ implies that the identifications in \eqref{transfer-T} are compatible with dualities)
we get the identifications
\be \lab{transfer-T}
(\Lie \cJ)_{\bZ_p} \cong H^1(X_0(n)_{\bZ_p}, \sO)  \qq \text{and} \qq  H^0(\cJ_{\bZ_p}, \Omega^1) \cong H^0(X_0(n)_{\bZ_p}, \Omega)
\ee
that are compatible with duality pairings, where $\Omega$ denotes the relative dualizing sheaf of $X_0(n)_{\bZ_p}$ over $\bZ_p$.
We transfer the $\bT$-action across these identifications to endow $H^1(X_0(n)_{\bZ_p}, \sO)$ and $H^0(\cJ_{\bZ_p}, \Omega^1)$ with a $\bT$-module structure.
\epp

The following result lies at the heart of our approach to the semistable case of the Manin conjecture.

\bprop \lab{H1-free}
For a maximal ideal $\fm \subset \bT$ of residue characteristic $p$, 
\benumr
\item \lab{H1f-i}
if $\ord_p(n) = 0$\uscolon or

\item \lab{H1f-iii}
if $\ord_p(n) = 1$ and $U_p \bmod \fm$ lies in $\bF_p^\times \subset \bT/\fm$ \up{the latter holds if $\fm$ contains the kernel of the map $q_f\colon \bT \surjects \cO_f$ determined by some newform $f$ of level $\Gamma_0(n)$  because $q_f(U_p) = \pm 1$}\uscolon or %\cite{AL78}*{p.~224 and Thm.~2.1}

\item \lab{H1f-ii}
if $\ord_p(n) = 1$ and $p$ is odd\uscolon
\eenum
then the following equivalent conditions hold\ucolon
\benuma
\item \lab{H1f-1}
the $\bT_\fm$-module $(\Lie \cJ)_{\bZ_p}\tensor_{\bT_{\bZ_p}} \bT_\fm$ is free of rank $1$\uscolon

\item \lab{H1f-2}
the $\bT_\fm$-module $H^1(X_0(n)_{\bZ_p}, \sO) \tensor_{\bT_{\bZ_p}} \bT_\fm$ is free of rank $1$\uscolon

\item \lab{H1f-3}
multiplicity one for differentials holds at $\fm$ in the sense that 
\[
\q \dim_{\bT/\fm} \p{H^0(X_0(n)_{\bF_p}, \Omega)[\fm]}  = 1.
\]
\eenum
\eprop

% Another approach to the proof of the last part was to work with dual of the sequence \eqref{Alb-seq}, but for $\Gamma_1(n)$ instead of $\Gamma_0(n)$ (for which these is an annoying problem of non-injectivity of $(\pi_\forg^*, \pi_\quot^*)$), and the to deduce the $\Gamma_0(n)$ case on the level of Hecke algebras using exactness on Neron models (and also that the dual of the Shimura group is constant). See the 11/10/16 write-up for an outline of this argument. This approach would have the advantage of proving the statement for any $\Gamma_H \cap \Gamma_0(p)$ level with $H$ intermediate between $\Gamma_1(\f{n}{p})$ and $\Gamma_0(\f{n}{p})$. The disadvantage is that one would have to first also prove the first two parts for such levels (but for this, Wiles' argument should extend without a problem; note also that in his argument there is no need to mess around with the minimal proper regular resolution because the $\cJ^0$ vs. $\Pic^0$ comparison also works for semistable curves (ARS seem to have missed this simplification)).

\bpf
The equivalence of \ref{H1f-1} and \ref{H1f-2} follows from \eqref{transfer-T}. By the formalism of cohomology and base change (see \cite{Ill05}*{8.3.11})  and Grothendieck--Serre duality (see \cite{Con00}*{Cor.~5.1.3}),
\[
H^1(X_0(n)_{\bZ_p}, \sO) \tensor_{\bZ_p} \bF_p \cong H^1(X_0(n)_{\bF_p}, \sO) \q \text{and} \q H^0(X_0(n)_{\bZ_p}, \Omega) \tensor_{\bZ_p} \bF_p \cong H^0(X_0(n)_{\bF_p}, \Omega);
\]
the latter identification endows $H^0(X_0(n)_{\bF_p}, \Omega)$ with the $\bT$-action used in \ref{H1f-3}.  Therefore,
\be \lab{m-dual}
H^1(X_0(n)_{\bZ_p}, \sO) \tensor_{\bT} \bT/\fm \qq \text{is the $\bF_p$-linear dual of} \qq H^0(X_0(n)_{\bF_p}, \Omega)[\fm],
\ee
and it follows that \ref{H1f-2} implies \ref{H1f-3}. Conversely, if \ref{H1f-3} holds, then, due to \eqref{m-dual}, the Nakayama lemma supplies a $\bT_\fm$-module surjection
\be \lab{s-Nak-arg}
s\colon \bT_\fm \surjects H^1(X_0(n)_{\bZ_p}, \sO) \tensor_{\bT_{\bZ_p}} \bT_\fm.
\ee
Since $H^1(X_0(n)_{\bZ_p}, \sO) \tensor_{\bT_{\bZ_p}} \bT_\fm$ is a faithful $\bT_\fm$-module (see \eqref{transfer-T}), the map $s$ is also injective, and hence is an isomorphism, which proves that \ref{H1f-3} implies \ref{H1f-2}.

The arguments above also apply to the minimal regular resolution $\wt{X_0(n)}_{\bZ_p}$ in place of $X_0(n)_{\bZ_p}$, so in the conditions \ref{H1f-2}--\ref{H1f-3} we could have instead used $\wt{X_0(n)}_{\bZ_p}$. Therefore, the results of \cite{ARS12}*{\S5.2} (which use $\wt{X_0(n)}_{\bZ_p}$), specifically, \cite{ARS12}*{Lemma 5.20},\footnote{The key inputs to the proof of loc.~cit.~are an Eichler--Shimura type congruence relation for $U_p$ in the style of \cite{Wil80}*{\S5} and arguments from \cite{Wil95}*{proof of Lemma 2.2} that use the $q$-expansion principle as in \cite{Maz77}*{pp.~94--95}.} show that either \ref{H1f-i} or \ref{H1f-iii} implies \ref{H1f-3}. Alternatively, \ref{H1f-i} implies \ref{H1f-1} by \cite{Par99}*{Thm.~4.2}.

The case \ref{H1f-ii} will only be used in \Cref{ARS-conj} and \Cref{Tpr}, so the cases \ref{H1f-i} and \ref{H1f-iii} suffice for the main results of the paper. To address the case \ref{H1f-ii}, we now assume that $p$ is odd with $\ord_p(n) = 1$, and we seek to show \ref{H1f-1}, that is, that $(\Lie \cJ)_\fm$ is free of rank $1$ as a $\bT_\fm$-module. Let 
\[
\tst \pi_\forg, \pi_\quot\colon X_0(n)_\bQ \rightrightarrows X_0(\f{n}{p})_\bQ
\]
be the degeneracy morphisms characterized as follows in terms of the moduli interpretation on the elliptic curve locus: $\pi_\forg$ forgets the $p$-primary factor of the cyclic subgroup of order $n$, whereas $\pi_\quot$ quotients the elliptic curve by this $p$-primary factor. We will consider the short exact sequence
\be \lab{Alb-seq}
\tst 0 \ra K \ra J_0(n) \xra{((\pi_\forg)_*,\, (\pi_\quot)_*)} J_0(\f{n}{p}) \times J_0(\f{n}{p}) \ra 0
\ee
in which $(\pi_\forg)_*$ and $(\pi_\quot)_*$ are induced by the Albanese functoriality, the surjectivity follows from \cite{Rib84}*{Cor.~4.2}, and $K$ is defined as the kernel. By loc.~cit.~and \cite{LO91}*{Thm.~2}, the component group scheme $K/K^0$ is constant, whereas the identity component $K^0$ is identified with the $p$-new subvariety of $J_0(n)$. The maps $(\pi_\forg)_*$ and $(\pi_\quot)_*$ commute with the Hecke operators $T_\ell$ and $U_\ell$ provided that $\ell \neq p$, % to see this draw the diagrams defining the Hecke correspondences and work away from the cusps; one of the squares in these diagrams (the left one) is then Cartesian on the level of $\ov{\bQ}$-points away from the cusps (or on the level of open modular stacks).
so they intertwine the actions of the ``$p$-anemic'' Hecke algebras 
\[
\tst \bT^{(p)} \subset \End(J_0(n)) \qq \text{and} \qq \bT_{p\text{-old}}^{(p)} \subset \End(J_0(\f{n}{p}))
\]
that are generated by these operators, and hence they define a surjective ring homomorphism
\be \lab{T-old-surj}
\bT^{(p)} \surjects \bT_{p\text{-old}}^{(p)}, \qq T_\ell \mapsto T_\ell, \q U_\ell \mapsto U_\ell.
\ee
Moreover, since $p$ is odd and does not divide $\f{n}{p}$, one knows from \cite{Wil95}*{Lemma on p.~491} that $T_p \in \bT_{p\text{-old}}^{(p)}$, that is, that $\bT_{p\text{-old}}^{(p)}$ is in fact the full Hecke algebra $\bT_{p\text{-old}} \subset \End(J_0(\f{n}{p}))$.

To determine the endomorphism of $J_0(\f{n}{p}) \times J_0(\f{n}{p})$ that intertwines $U_p \in \End(J_0(n))$, we regard $Y_0(n)_\bQ$ as the coarse moduli space of pairs 
\[
(\phi\colon E_1 \ra E_2, C \subset E_1),
\]
where $\phi$ is a $p$-isogeny of elliptic curves and $C$ is a cyclic subgroup of order $\f{n}{p}$, so that the $p$-Atkin--Lehner involution $w_p$ sends $(\phi, C)$ to $(\phi^\vee, \phi(C))$. In terms of this interpretation, $U_p$ quotients $E_1$ by a variable subgroup $C_p' \subset E_1$ of order $p$ such that $C_p' \cap \Ker \phi = 0$. Therefore, $U_p + w_p$ sends
\[
(\phi\colon E_1 \ra E_2, C \subset E_1) \qq \text{to} \qq \sum (\psi\colon E_3 \ra E_1, \psi\i(C) \subset E_3),
\]
where the sum runs over the isomorphism classes of all $p$-isogenies that cover $E_1$. In particular,\footnote{Compare with the formula in \cite{Rib90}*{proof of Prop. 3.7} that used the ``Picard functoriality'' Hecke operators.}
\[
U_p + w_p = (\pi_\quot)^* \circ (\pi_\forg)_* \qq \text{inside} \qq \End(J_0(n)).
\]
Thus, since $T_p = (\pi_\forg)_* \circ (\pi_\quot)^*$ in $\End(J_0(\f{n}{p}))$ (the switch to ``Picard functoriality'' here does not matter because $(p, \f{n}{p}) = 1$, see \cite{MW84}*{Ch.~II, \S5.4 (2)}) and $ (\pi_\quot)_* \circ (\pi_\quot)^* = p + 1$,
\be \lab{Up-inter}
\tst U_p \q \text{intertwines the endomorphism} \q (x, y) \mapsto (T_p x - y, px) \q \text{of} \q J_0(\f{n}{p}) \times J_0(\f{n}{p}).
\ee
% Compare with \cite{Wil95}*{p.~490}.
The $p$-new subvariety $K^0$ is isogenous to a product of newform quotients (with multiplicities) of variable $J_0(n')$ for divisors $n'\mid n$ such that $n' \nmid \f{n}{p}$, that is, such that $p \mid n'$. Since the $U_p$ operator commutes with the degeneracy maps towards such $J_0(n')$, it acts as $\pm 1$ on each simple isogeny factor of $K^0$. In particular, since $(\Lie J_0(n))_{\bQ_p} \simeq \bT_{\bQ_p}$ as $\bT_{\bQ_p}$-modules (see \cite{DDT97}*{1.34}), the factors $\bT_{\fm'}$ of $\bT_{\bZ_p}$ that meet the support of $(\Lie K)_{\bQ_p}$ must satisfy either $U_p + 1 \in \fm'$ or $U_p - 1 \in \fm'$. By \ref{H1f-iii}, we may assume $\bT_\fm$ is not such factor, to the effect that 
\be \lab{m-seems-old}
(\Lie K) \tensor_{\bT_{\bZ_p}} \bT_\fm = 0.
\ee
Since $p$ is odd and $K/K^0$ is constant, the maps $J_0(n) \ra J_0(n)/K^0$ and $J_0(n)/K^0 \ra J_0(n)/K$ induce smooth morphisms on N\'{e}ron models over $\bZ_p$ (see \cite{BLR90}*{7.5/4 (ii) and its proof, 7.5/6}). Thus, the sequence \eqref{Alb-seq} induces a short exact sequence on Lie algebras of N\'{e}ron models over $\bZ_p$:
\[
0 \ra (\Lie \cK)_{\bZ_p} \ra (\Lie \cJ)_{\bZ_p} \ra (\Lie \cJ' \times \Lie \cJ')_{\bZ_p} \ra 0.
\]
The vanishing \eqref{m-seems-old} then implies that 
\be \lab{loc-m-iso}
(\Lie \cJ)_\fm \isomto (\Lie \cJ' \times \Lie \cJ')_\fm,
\ee
so the maximal ideal $\fm^{(p)} \ce \fm \cap \bT^{(p)}$ of $\bT^{(p)}$ is such that $(\Lie \cJ' \times \Lie \cJ')_{\fm^{(p)}} \neq 0$, to the effect that $\fm^{(p)}$ is ``$p$-old,'' that is, is identified with a maximal ideal $\fm^{(p)} \subset \bT_{p\text{-old}}$ via \eqref{T-old-surj} (and the sentence that follows \eqref{T-old-surj}). Thus, the case \ref{H1f-i} implies that the $\bT^{(p)}/\fm^{(p)}$-vector space $(\Lie \cJ' \times \Lie \cJ')/\fm^{(p)}$ is of dimension $2$. Moreover, due to the formula \eqref{Up-inter}, the $U_p$ operator does not act as a $\bT^{(p)}/\fm^{(p)}$-scalar on this vector space, so the $\bT/\fm$-vector space $(\Lie \cJ)/\fm$, which, due to \eqref{loc-m-iso}, is a further quotient of $(\Lie \cJ' \times \Lie \cJ')/\fm^{(p)}$, is of dimension $\le 1$. It then follows from the Nakayama lemma as in the paragraph after \eqref{s-Nak-arg} that $(\Lie \cJ)_\fm$ is free of rank $1$ as a $\bT_\fm$-module, as desired.
\epf

\brem \lab{ARS-conj}
In \cite{ARS12}*{\S5.2.1}, based on computational evidence, Agashe, Ribet, and Stein asked whether \ref{H1f-ii} implies \ref{H1f-3} and showed that this is not the case if the parity condition is dropped in \ref{H1f-ii}. There they also implicitly raised the question answered by part \ref{Tpr-ii} of the following corollary.
\erem

\bcor\lab{Tpr}
For a maximal ideal $\fm \subset \bT$ of residue characteristic $p$,
\benumr
\item \lab{Tpr-i}
if $\ord_p(n) = 0$\uscolon or

\item \lab{Tpr-iii}
if $\ord_p(n) = 1$ and $U_p \bmod \fm$ lies in $\bF_p^\times \subset \bT/\fm$\uscolon or

\item \lab{Tpr-ii}
if $\ord_p(n) = 1$ and $p$ is odd\uscolon
\eenum
then the saturation $\bT' \ce \bT_\bQ \cap \End_\bQ(J_0(n))$ of $\bT$ agrees with $\bT$ at $\fm$, that is, 
\[
\xymatrix{\bT_\fm \ar@{^(->}[r]^{\sim} &\bT'_\fm.}
\] 
In particular, if $n$ is odd and squarefree, then the inclusion $\bT \hra \bT'$ is an isomorphism.
\ecor

\bpf
Since $\bT'$ acts on $\Lie \cJ$ faithfully and $\bT$-linearly, $\bT\pr_\fm$ acts on 
\[
(\Lie \cJ)_{\bZ_p}\tensor_{\bT_{\bZ_p}} \bT_\fm\overset{\ref{H1-free}}{\simeq} \bT_\fm
\]
faithfully and $\bT_\fm$-linearly, and the desired conclusion follows.
\epf

\brem
When $n$ is a prime, the last claim of \Cref{Tpr} also follows from \cite{Maz77}*{II.9.5}.
\erem

For our purposes, the significance of \Cref{H1-free}, especially of its case \ref{H1f-iii}, is the resulting control of the congruence module $\mathrm{cong}_{\Lie \cJ}$ introduced in \Cref{cong-mod-M} with the following setup.

\bpp[A modular elliptic curve] \lab{mod-ell-curve}
For the rest of \S\ref{ell-curve} we
\begin{itemize}
\item
fix a new elliptic optimal quotient $\pi\colon J_0(n) \surjects E$;

\item
let $f \in H^0(X_0(n)_{\bQ}, \Omega^1)$ be the normalized newform determined by $\pi$;

\item
let $e_f \in \bT_\bQ$ be the idempotent that cuts out the factor of $\bT_\bQ$ that corresponds to $f$;

\item
let $e_{f^\perp} \ce 1 - e_f$ be the complementary idempotent in $\bT_\bQ$.
\end{itemize}

The idempotents $e_f$ and $e_{f^\perp}$ decompose every $\bT$-module $M$ rationally: 
\[
M_\bQ \cong M_\bQ[e_f] \oplus M_\bQ[e_{f^\perp}].
\]
The following congruence module measures the failure of an analogous integral decomposition.
\epp

\bd \lab{cong-mod-M}
The \emph{congruence module} of a $\bZ$-torsion free $\bT$-module $M$ is the quotient (see \eqref{abuse})
\[
\tst    \mathrm{cong}_M \ce \f{M}{M[e_f] + M[e_{f^\perp}]}.
\]
\ed

\beg \lab{cong-mod-eg}
With the choice $M = \prod_p H^1_\et(J_0(n)_{\ov{\bQ}}, \bZ_p)$, the congruence module $\mathrm{cong}_M$ is isomorphic to $(\bZ/\deg_f \bZ)^2$, where the \emph{modular degree} $\deg_f$ is the positive integer that equals $\pi \circ \pi^\vee$ in $\End_\bQ(E)$ (so that $e_f = \f{\pi \circ \pi^\vee}{\deg_f}$). Indeed, this follows from the optimality of $\pi$ and from the observation that $H^1_\et((-)_{\ov{\bQ}}, \bZ_p)$ carries short exact sequences of abelian varieties to those of finite free $\bZ_p$-modules.
\eeg

This example leads to the following lemma, whose proof is a variant of the proof of \cite{DDT97}*{Lem.~4.17} (the quotient $\cO/\eta_\Sigma$ used there is a congruence module). The lemma is well known and also follows from \cite{AU96}*{Lem.~3.2}, \cite{CK04}*{Thm.~1.1}, or \cite{ARS12}*{Thm.~2.1}. 

\blem \lab{deg-cong-div}
In the setup of \uS\uref{mod-ell-curve}, we have $\deg_f \mid \#(\mathrm{cong}_{\bT})$.
\elem

\bpf
For every prime $p$, the $\bT_{\bQ_p}$-module $H^1_\et(J_0(n)_{\ov{\bQ}}, \bQ_p)$ is free of rank $2$ (see \cite{DDT97}*{Lem.~1.38--1.39}), so the module 
\[
\tst \f{H^1_\et(J_0(n)_{\ov{\bQ}},\, \bZ_p)}{H^1_\et(J_0(n)_{\ov{\bQ}},\, \bZ_p)[e_f]} \qq \text{is free of rank $2$ over} \qq (\bT/\bT[e_{f}])_{\bZ_p} \cong \bZ_p.
\]
In particular, $\mathrm{cong}_{\bT_{\bZ_p}}^{\oplus 2}$ surjects onto $\mathrm{cong}_{H^1_\et(J_0(n)_{\ov{\bQ}},\, \bZ_p)}$ and the claim follows from \Cref{cong-mod-eg}.
\epf

We are ready to exploit \Cref{H1-free} in the proof of the following key lemma, which will give us the semistable case of the Manin conjecture in \Cref{Manin-sst} and whose proof will simultaneously reprove \cite{ARS12}*{Thm.~2.1}. An alternative route would be to deduce \Cref{Manin-sst} from the results of \S\ref{contorni}, e.g.,~from \Cref{genl-results}, that are valid for newform quotients of arbitrary dimension.

\blem \lab{cong-deg}
For every prime $p$ such that $p^2 \nmid n$,
\be \lab{CD-eq}
\ord_p(\# (\mathrm{cong}_{\Lie \cJ})) = \ord_p(\deg_f),
\ee
the map $(\Lie \cJ)_{\bZ_p} \ra (\Lie \cE)_{\bZ_p}$ is surjective, and $\xymatrix{((\Lie \pi^\vee)(\Lie \cE))_{\bZ_p} \ar@{^(->}[r]^-{\sim} & (\Lie \cJ)_{\bZ_p}[e_{f^\perp}].}$
\elem

\bpf
The formation of the $\bT$-module 
\[
\tst \f{\Lie \cJ}{(\Lie \cJ)[e_f] + (\Lie \cJ)[e_{f^\perp}]}
\]
commutes with (flat) base change to $\bZ_p$, so we let $\fn$ range over the maximal ideals of $\bT$ of residue characteristic $p$ and decompose $\bT_{\bZ_p} \cong \prod \bT_{\fn}$. 

Let $\fm \subset \bT$ be the preimage of $p\bZ$ under the surjection $\bT \surjects \bZ$ determined by $f$. The image of $e_f$ in $\bT_{\bQ_p}$ lies in $(\bT_\fm)_{\bQ_p}$, so for every $\fn \neq \fm$ the ``$(\bT_\fn)_{\bQ_p}$-coordinate'' of $e_{f}$ vanishes, to the effect that 
\[
\bT_\fn = \bT_\fn[e_f] \qq \text{and} \qq (\Lie \cJ) \tensor_{\bT_{\bZ_p}} \bT_{\fn} = ((\Lie \cJ) \tensor_{\bT_{\bZ_p}} \bT_{\fn})[e_f]
\]
for every such $\fn$. Consequently, 
\[
\tst 
\p{\f{\Lie \cJ}{(\Lie \cJ)[e_f] + (\Lie \cJ)[e_{f^\perp}]} } \tensor_\bZ \bZ_p \cong \f{(\Lie \cJ)_{\fm}}{(\Lie \cJ)_{\fm}[e_f] + (\Lie \cJ)_{\fm}[e_{f^\perp}]} \overset{\ref{H1-free}}{\simeq} \f{\bT_\fm}{\bT_\fm[e_f] + \bT_\fm[e_{f^\perp}]} \cong \p{\f{\bT}{\bT[e_f] + \bT[e_{f^\perp}]}} \tensor_\bZ \bZ_p,
\]
so that  
\[
\ord_p(\#(\mathrm{cong}_{\Lie \cJ})) = \ord_p(\#(\mathrm{cong}_\bT)).
\]
When combined with \Cref{deg-cong-div}, this gives the inequality ``$\ge$'' in \eqref{CD-eq}.

For the converse inequality, we let $\cE$ denote the N\'{e}ron model of $E$ over $\bZ$, observe that the injection
\[
\tst \Lie \pi \colon \f{\Lie \cJ}{(\Lie \cJ)[e_f]} \hra \Lie \cE \qq \text{gives} \qq \f{\Lie \cJ}{(\Lie \cJ)[e_{f}] + (\Lie \pi^\vee)(\Lie \cE)} \hra \f {\Lie\cE}{\deg_f \cdot \Lie \cE} \simeq \f{\bZ}{\deg_f \cdot \bZ},
\]
and conclude by using the inclusion $(\Lie \pi^\vee)(\Lie \cE) \subset (\Lie \cJ)[e_{f^\perp}]$.
\epf

\bthm \lab{Manin-sst}
For a new elliptic optimal quotient $\pi\colon J_0(n) \surjects E$, the Manin constant $c_\pi$ satisfies
\[
\ord_p(c_\pi) = 0 \qq \text{for every prime} \q p \q \text{such that} \q p^2 \nmid n. 
\]
\ethm

\bpf
By \Cref{cong-deg}, the map $(\Lie \cJ)_{\bZ_p} \ra (\Lie \cE)_{\bZ_p}$ is surjective. Thus, $\f{H^0(\cJ_{\bZ_p}, \Omega^1)}{\pi^*(H^0(\cE_{\bZ_p}, \Omega^1))}$ is torsion free, that is, that the pullback $\pi^*(\omega_E)$ of a N\'{e}ron differential of $E$ is not divisible by $p$ in $H^0(\cJ_{\bZ_p}, \Omega^1)$.

The $p$-Atkin--Lehner involution sends $f$ to $\pm f$, so $f$ lies in 
\[
H^0(X_0(n)_{\bZ_p}, \Omega) \overset{\eqref{transfer-T}}{\cong} H^0(\cJ_{\bZ_p}, \Omega^1)
\]
(see \cite{Ces16g}*{2.7--2.8}). Moreover, by definition, $\pi^*(\omega_E) = c_\pi \cdot f$ and, by \cite{Edi91}*{Prop.~2}, we have $c_\pi \in \bZ$. Thus, since $\pi^*(\omega_E)$ is not divisible by $p$ in $H^0(\cJ_{\bZ_p}, \Omega^1)$, the desired $\ord_p(c_\pi) = 0$ follows.
\epf

We end \S\ref{ell-curve} with the semistable case of the analogue of the Manin conjecture for parametrizations by modular curves intermediate between $X_1(n)$ and $X_0(n)$ (see \Cref{manin-H-sst}). The following variant of \cite{Ces16g}*{Lem.~4.4} (or of \cite{GV00}*{Prop.~3.3}) reduces this analogue to the Manin conjecture itself.

\blemt \lab{down}
Let $H, H' \subset \GL_2(\wh{\bZ})$ be subgroups such that $\Gamma_1(n) \subset H \subset H' \subset \Gamma_0(n)$ for some $n \in \bZ_{\ge 1}$, and let 
\[
\pi \colon J_H \surjects E \qq \text{and} \qq \pi' \colon J_{H'} \surjects E'
\]
be new elliptic optimal quotients such that $E$ and $E'$ are isogenous over $\bQ$. There is a unique isogeny $e$ making the diagram
\be \ba \lab{down-sq}
\xymatrix{
J_{H} \ar@{->>}[r]^-{\pi} \ar@{->>}[d]^-{j^\vee} & E \ar@{->>}[d]^-{e} \\
J_{H'} \ar@{->>}[r]^-{\pi'} & E'
}
\ea \ee
commute, where $j^\vee$ is the dual of the pullback map $j \colon J_{H'} \ra J_{H}$. The kernel $\Ker e$ is constant and is a subquotient of the Cartier dual of the Shimura subgroup $\Sigma(n) \subset J_0(n)$.  The Manin constants $c_{\pi}$ and $c_{\pi'}$ are nonzero integers related by the equality 
\be\lab{go-down}
c_{\pi'} = c_{\pi} \cdot \#\Coker\p{\Lie e\colon \Lie \cE \ra \Lie \cE'},
\ee
where $\cE$ and $\cE'$ are the N\'{e}ron models over $\bZ$ of $E$ and $E'$, respectively.
\elemt

\bpf
The existence (resp.,~uniqueness) of $e$ follows from the multiplicity one theorem (resp.,~from the surjectivity of $\pi$). For the rest, we loose no generality by assuming that $H' = \Gamma_0(n)$. Then $\Ker(e)^\vee$ is a subgroup of 
\[
\Sigma(n) \ce \Ker(J_0(n) \ra J_1(n)).
\]
Since $\Sigma(n)$ is of multiplicative type (see \cite{LO91}*{Thm.~2}), the claims about $\Ker e$ follow. The formation of the normalized newform  $f \in H^0(J_{H'}, \Omega^1)$ determined by the isogeny class of $E$ is compatible with pullback by $j^\vee$ (see \cite{Ces16g}*{2.8}), so the comparison of the two ways to pull back a N\'{e}ron differential of $E'$ to $J_{H}$ in \eqref{down-sq} gives \eqref{go-down}. For the integrality of $c_\pi$ it then suffices to assume that  $H = \Gamma_1(n)$ and to apply \cite{Ste89}*{Thm.~1.6}. 
\epf

\bthm \lab{manin-H-sst}
For an $n \in \bZ_{\ge 1}$, a subgroup $H \le \GL_2(\wh{\bZ})$ such that $\Gamma_1(n) \subset H \subset \Gamma_0(n)$, and a new elliptic optimal quotient $\pi \colon J_H \surjects E$, if $p$ is a prime with $p^2 \nmid n$, then the Manin constant $c_\pi$ satisfies $\ord_p(c_\pi) = 0$ and $\pi$ induces a smooth morphism $(\cJ_H)_{\bZ_p} \ra \cE_{\bZ_p}$ between the N\'{e}ron models over $\bZ_p$.
\ethm

\bpf
Let $\pi' \colon J_0(n) \surjects E'$ be the new elliptic optimal quotient for which $E'$ is isogenous to $E$ and let $e \colon E \ra E'$ be the isogeny supplied by \Cref{down}. By \Cref{Manin-sst}, $\ord_p(c_{\pi'}) = 0$, so, since $c_\pi \in \bZ$, \eqref{go-down} implies that $\ord_p(c_\pi) = 0$. Thus, since the elements of $H^0(\cJ_H, \Omega^1)$ have integral $q$-expansions (see \cite{CES03}*{proof of Lem.~6.1.6}), it follows that the pullback $\pi^*(\omega_E)$ of a N\'{e}ron differential of $E$ is not divisible by $p$ in $H^0((\cJ_H)_{\bZ_p}, \Omega^1)$. Consequently, $(\Lie \cJ_H)_{\bZ_p} \ra (\Lie \cE)_{\bZ_p}$ is surjective and the smoothness claim follows.
\epf

One consequence of \Cref{manin-H-sst} is the following additivity of Faltings height.

\bcort \lab{Faltings-ht}
For an $H$ as in Theorem \uref{manin-H-sst} and a new elliptic optimal quotient $\pi \colon J_H \surjects E$ such that $E$ is semistable, the Faltings height $h(-)$ over $\bQ$ satisfies 
\[
h(J_H) =  h(\Ker \pi) + h(E).
\]
\ecort

\bpf
\Cref{manin-H-sst} applies at every $p$, so it implies that the map $\cJ_H \ra \cE$ is smooth and hence, by \cite{BLR90}*{7.1/6}, that the N\'{e}ron model $\cK$ over $\bZ$ of $\Ker \pi$ is its kernel. It follows that the sequences
\[
0 \ra \Lie \cK \ra \Lie \cJ_H \ra \Lie \cE \ra 0 \qq \text{and} \qq 0 \ra H^0(\cE, \Omega^1) \ra H^0(\cJ_H, \Omega^1) \ra H^0(\cK, \Omega^1) \ra 0
\]
are short exact, so the arguments from \cite{Ull00}*{proof of Prop.~3.3} give the conclusion.
\epf

%%%%%%%%%%%%%%%%%%%%%%%%%%%%%%%%%%%%%%%

%\appendix

\section{A link between the integral $p$-adic \'{e}tale and de Rham cohomologies} \lab{et-dR}

We will analyze the semistable case of a generalization of the Manin conjecture to higher dimensional newform quotients in \S\ref{contorni}, and this will rest on the present section which establishes relations between the integral $p$-adic \'{e}tale and de Rham cohomologies of abelian varieties over $p$-adic fields. Our point of view is that arithmetic duality in the form of a result of Raynaud \cite{Ray85}*{Thm.~2.1.1}, which we recast and mildly generalize in \Cref{Raynaud-re,Raynaud-re-re}, is capable of supplying such relations.

%work of the present section whose goal is to relate exactness properties on  of a complex of abelian varieties over a $p$-adic field (the complexes that appear most frequently are isogenies and short exact sequences, but see also \Cref{BMS-eg} and its proof). 

\bpp[The field $K$] \lab{et-dR-not}
Throughout \S\ref{et-dR} we fix a mixed characteristic $(0, p)$ complete discretely valued field $K$ whose residue field $k$ is perfect. We denote the ring of integers of $K$ by $\cO$ and we denote N\'{e}ron models over $\cO$ by calligraphic letters: for instance, $\cA$ and $\cB$ are the N\'{e}ron models over $\cO$ of abelian $K$-varieties $A$ and $B$, whereas $\cA^\vee$ is the N\'{e}ron model of the dual abelian variety $A^\vee$. 
\epp

The following de Rham lattice $H^1_\dR(-/\cO)$ constructed by Mazur and Messing will be key for our purposes. In order to emphasize its functoriality, we review its definition.

\bpp[The integral $H^1_\dR$ of an abelian variety] \lab{MM-review}
Let $A$ be an abelian variety over $K$. By \cite{MM74}*{I.5.1},  the N\'{e}ron model $\cA^\vee$ of the dual abelian variety $A^\vee$ is identified with the fppf sheaf
\[
\sE xt^1(\cA^0, \bG_m) \qq\text{defined by} \qq  S \mapsto \Ext^1_S(\cA^0_S, (\bG_m)_S)
\]
on the category of smooth $\cO$-schemes $S$, where $\Ext^1_S(\cA^0_S, (\bG_m)_S)$ is the abelian group of extensions $0 \ra (\bG_m)_S \ra \cE \ra \cA^0_S \ra 0$ of fppf sheaves of abelian groups.\footnote{For the sheaf condition of $\sE xt^1(\cA^0, \bG_m)$, it is key that every $0 \ra (\bG_m)_S \ra \cE \ra \cA^0_S \ra 0$ has no nonidentity automorphisms, as may be checked over $S_K$ due to the $\cO$-flatness of $S$ (each $\cE$ is necessarily a smooth $S$-scheme).} One also considers the sheaf 
\[
\sE xtrig^1(\cA^0, \bG_ m) \qq\text{defined by} \qq  S \mapsto \Extrig^1_S(\cA^0_S, (\bG_m)_S)
\]
on the category of smooth $\cO$-schemes $S$, where $\Extrig^1_S(\cA^0_S, (\bG_m)_S)$ is the abelian group of extensions as before that are, in addition, equipped with an $S$-morphism $r$, a rigidification, that fits into a commutative diagram
\[
\xymatrix{
 & &  & \underline{\Spec}(\sO_{\cA^0_S}/\sI^2_{e_{\cA^0_S}}) \ar@{^(->}[d]_-{i} \ar[ld]_-{r} & \\
0 \ar[r] & (\bG_m)_S \ar[r] & \cE \ar[r] & \cA_S^0 \ar[r] & 0
}
\]
in which $i$ is the first infinitesimal neighborhood of the identity section of $\cA_S^0$. By \cite{MM74}*{I.5.2}, the functor $\sE xtrig^1(\cA^0, \bG_ m)$ is representable by a smooth $\cO$-group scheme that fits into a short exact sequence
\be \lab{extrig-seq}
0 \ra \bV(\Lie \cA) \ra \sE xtrig^1(\cA^0, \bG_ m) \ra \cA^\vee \ra 0
\ee
in which $\bV(\Lie \cA)$ is the affine smooth $\cO$-group scheme that represents the finite free $\cO$-module $(\Lie \cA)^*$ (see \cite{SGA3Inew}*{I, 4.6.3.1}). By \cite{MM74}*{I.2.6.7}, the $K$-base change of the sequence \eqref{extrig-seq} is identified with the universal vector extension of $A^\vee$, so 
\[
\Lie(\sE xtrig^1(\cA^0, \bG_ m)_K) \qq \text{is identified with} \qq H^1_\dR(A/K).
\]
Therefore, the finite free $\cO$-module 
\[
H^1_\dR(\cA/\cO) \ce \Lie(\sE xtrig^1(\cA^0, \bG_ m))
\]
is a natural integral structure on $H^1_\dR(A/K)$. If $A$ has good reduction, that is, if $\cA$ is an abelian scheme, then, by \cite{MM74}*{I.\S4, esp.,~I.3.2.3 a), I.4.2.1, and I.4.1.7}, the $\cO$-module $H^1_\dR(\cA/\cO)$ is identified with the first de Rham cohomology group of $\cA$ over $\cO$. By construction, the formation of \eqref{extrig-seq} and of $H^1_\dR(\cA/\cO)$ is contravariantly functorial in $A$ and so is the formation of the exact sequence
\be \lab{filt}
0 \ra (\Lie \cA)^* \ra H^1_\dR(\cA/\cO) \ra \Lie \cA^\vee \ra 0
\ee
of finite free $\cO$-modules obtained from \eqref{extrig-seq}. 
\epp

\bpp[The normalized length $\val_\cO(-)$] \lab{norm-length}
For a finitely generated torsion $\cO$-module $M$, we set
\[
\tst \val_\cO(M) \ce \f{1}{e(K/\bQ_p)} \cdot \length_\cO(M), \qq \text{where} \qq e(K/\bQ_p) \ce \length_{\cO}(\cO/p)
\]
is the absolute ramification index of $K$. For a bounded complex $(M_\bullet, d_\bullet)$ of finitely generated $\cO$-modules such that $(M_\bullet)_K$ is exact, we set
\be \lab{norm-length-cx}
\tst \val_\cO(M_\bullet) = \sum_i (-1)^i \cdot \val_{\cO}\p{\f{\Ker(d_{i} \colon M_i \ra M_{i - 1})}{\im(d_{i + 1} \colon M_{i + 1} \ra M_{i})}}.
\ee
We prefer the formalism of normalized length $\val_\cO(-)$ to that of length $\length_\cO(-)$ because the former is insensitive  to base change to the ring of integers of a finite extension of $K$.
\epp

We are ready for the following variant of \cite{Ray85}*{Thm.~2.1.1} that (through the proof of loc.~cit.)~uses arithmetic duality results of B\'{e}gueri \cite{Beg80}. In its formulation, with an eye towards applications to modular Jacobians, we keep track of an action of a ``Hecke algebra'' $\bT$; the basic case is $\bT = \bZ$.

\bthm\lab{Raynaud-re}
Let $\bT$ be a commutative ring that is finite free as a $\bZ$-module, let $A$ and $B$ be abelian varieties over $K$ endowed with a $\bT$-action, and let $f \colon A \ra B$ be a $\bT$-equivariant $K$-isogeny. For every maximal ideal $\fm \subset \bT$ of residue characteristic $p$, we have
\be \lab{RR-eq}
\tst \val_{\bZ_p} \p{\f{H^1_\et(A_{\ov{K}},\, \bZ_p)_\fm}{f^*(H^1_\et(B_{\ov{K}},\, \bZ_p)_\fm)}} =  \val_\cO\p{\f{H^1_\dR(\cA/\cO)_\fm}{f^*(H^1_\dR(\cB/\cO)_\fm)}};
\ee
in addition, both sides of this equality are equal to $\ord_p(\#(\Ker f)[\fm^\infty])$.
\ethm

\bpf
Since $H^1_\et(A_{\ov{K}}, \bZ_p)$ is identified with the $\bZ_p$-linear dual of the $p$-adic Tate module $T_p(A_{\ov{K}})$ compatibly with the $\bT_{\bZ_p}$-action, and likewise for $B$, the left side of \eqref{RR-eq} equals $\ord_p(\#(\Ker f)[\fm^\infty])$.

If we ignore the $\bT$-action, more precisely, if we take $\bT = \bZ$, then \cite{Ray85}*{Thm.~2.1.1} gives
\be \lab{Ray-input}
\tst  \val_\cO\p{\f{\Lie \cB}{(\Lie f)(\Lie \cA)}} + \val_\cO\p{\f{\Lie \cA^\vee}{(\Lie f^\vee)(\Lie \cB^\vee)}} = \val_\cO(\f{\cO}{(\#(\Ker f))}) = \ord_p(\#(\Ker f)).
\ee
Thus, in the $\bT = \bZ$ case \eqref{RR-eq}  follows from \eqref{Ray-input} and from the commutative diagram
\be\ba \lab{RR-diag}
\xymatrix{
0 \ar[r] & (\Lie \cB)^* \ar[r] \ar@{^(->}[d]^{(\Lie f)^*} & H^1_\dR(\cB/\cO) \ar[r] \ar@{^(->}[d]^{f^*} & \Lie \cB^\vee \ar[r]\ar@{^(->}[d]^{\Lie(f^\vee)} & 0 \\
0 \ar[r] & (\Lie \cA)^* \ar[r] & H^1_\dR(\cA/\cO) \ar[r] & \Lie \cA^\vee \ar[r] & 0.
}
\ea \ee
In the general case, both sides of \eqref{RR-eq} are additive in composites of isogenies, so we factor $f$ to assume that $\Ker f$, and hence also the left side of \eqref{RR-eq}, is supported entirely at $\fm$. Then we use the decomposition $\bT_{\bZ_p} \cong \prod_{\fn} \bT_{\fn}$, where $\fn$ ranges over the maximal ideals of $\bT$ of residue characteristic $p$, to find a $t \in \bT$ that kills $\Ker f$ but pulls back to a unit in every $\bT_{\fn}$ with $\fn \neq \fm$ and to a unit in $\bT_{\bZ_\ell}$ for some fixed auxiliary prime $\ell \neq p$. This $\ell$-adic assumption ensures that multiplication by $t$ is a self-isogeny of $A$ (of degree prime to $\ell$), whereas the inclusion $\Ker f \subset \Ker t$ translates into a factorization $t = g \circ f$ for some isogeny $g \colon B \ra A$. It then follows from the choice of $t$ that for every $\fn \neq \fm$ the injection 
\[
f^*\colon H^1_\dR(\cB/\cO)_{\fn} \hra H^1_\dR(\cA/\cO)_{\fn}
\]
must be surjective. In conclusion, the quotient $\f{H^1_\dR(\cA/\cO)}{f^*(H^1_\dR(\cB/\cO))}$ is also supported entirely at $\fm$, to the effect that \eqref{RR-eq} follows from its special $\bT = \bZ$ case.
\epf

\brem
Both sides of \eqref{RR-eq} are additive in composites of isogenies and are equal for the multiplication by $n$, so ``$\le$'' in \eqref{RR-eq} implies the equality. In the good reduction case this inequality may be deduced from integral $p$-adic Hodge theory (instead of arithmetic duality \cite{Ray85}*{Thm.~2.1.1}): one may use the results of \cite{BMS16} to build an $A_\Inf$-module $\f{H^1_{A_\Inf}(\cA)}{f^*(H^1_{A_\Inf}(\cB))}$ that in the key $\bT = \bZ$ case realizes the quotient on the right side of \eqref{RR-eq} as a specialization of that on the left side.
\erem

In order to make \Cref{Raynaud-re} applicable more broadly, we turn to its following variant.

\bthm \lab{Raynaud-re-re}
Let 
\[
A_\bullet = \ldots \ra A_{i - 1} \xra{d_{i-1}} A_i \xra{d_i} A_{i + 1} \ra \ldots,
\]
where $d_i \circ d_{i -1} = 0$ for every $i$, be a complex of abelian varieties over $K$. Suppose that
\benumr
\item \lab{RR-i}
the complex $A_\bullet$ is \emph{bounded} in the sense that $A_i = 0$ for all but finitely many $i$\uscolon and

\item \lab{RR-ii}
the complex $A_\bullet$ is \emph{exact up to isogeny} in the sense that $\im d_{i - 1} = (\Ker d_i)^0$ for every $i$.
\eenum
Under these assumptions, with the notation of \eqref{norm-length-cx} we have
\be \lab{RRR-eq}
\tst \val_{\bZ_p}\p{H^1_\et((A_\bullet)_{\ov{K}},\, \bZ_p)} = \val_{\cO}\p{H^1_\dR(\cA_\bullet/\cO)},
\ee
where the complex $H^1_\et((A_\bullet)_{\ov{K}},\, \bZ_p)$ has the term $H^1_\et((A_i)_{\ov{K}},\, \bZ_p)$ in degree $i$ and similarly for $H^1_\dR(\cA_\bullet/\cO)$\uscolon in addition, both sides of \eqref{RRR-eq} are equal to
\[
\tst \sum_i (-1)^i \ord_p\p{\#\p{\f{\Ker d_i}{\im d_{i - 1}}}},
\]
so that, in particular, if $A_\bullet$ is exact, then both sides of \eqref{RRR-eq} vanish.
\ethm

The reduction of \Cref{Raynaud-re-re} to \Cref{Raynaud-re} rests on the following lemma, which, in addition, shows that after inverting $p$ the \'{e}tale and the de Rham complexes appearing in \eqref{RRR-eq} are exact.

\blem \lab{really-isog} 
For an $A_\bullet$ that satisfies \ref{RR-i} and \ref{RR-ii}, there is a $K$-morphism 
\[
f_\bullet \colon (\wt{A}_\bullet, \wt{d}_\bullet) \ra (A_\bullet, d_\bullet)
\]
such that each $f_i$ is an isogeny and $\wt{A}_\bullet$ is a \emph{split} complex of abelian varieties over $K$ in the sense that there are $K$-isomorphisms $\wt{A}_i \cong \Ker \wt{d}_i \times_K \Ker \wt{d}_{i + 1}$ that are compatible with the differentials $\wt{d}_i$.
\elem

\bpf
The maps $d_i \colon A_i \ra \im d_i$ are surjective, so the Poincar\'{e} complete reducibility theorem (see \cite{Con06}*{Cor.~3.20}) supplies an abelian subvariety $B_i \subset A_i$ for which $(d_{i})|_{B_i} \colon B_i \ra \im d_i$ is an isogeny. By letting $f_i$ be the sum of $(d_{i - 1})|_{B_{i - 1}}$ and of $B_i \hra A_i$, we obtain the commutative diagram
\[
\xymatrix{
\dots \ar[r] & B_{i - 1} \times_K B_i \ar[d]^{f_i} \ar[r]^-{\wt{d}_i} & B_{i} \times_K B_{i + 1} \ar[r] \ar[d]^{f_{i + 1}} & \dots \\
\dots \ar[r] & A_i \ar[r]^{d_i} & A_{i + 1} \ar[r] & \dots
}
\]
whose top row is the candidate $\wt{A}_\bullet$. Each $\Ker f_i$ is finite because so is each $\im d_{i - 1} \cap B_i$, and for every $i$ we have 
\[
\dim A_i = \dim B_i + \dim B_{i - 1}.
\]
In conclusion, each $f_i$ is an isogeny. 
\epf

\emph{Proof of Theorem \uref{Raynaud-re-re}.} 
\Cref{really-isog} supplies a split bounded complex $\wt{A}_\bullet$ of abelian varieties over $K$ and a $K$-isogeny $f_\bullet \colon \wt{A}_\bullet \ra A_\bullet$. For this split $\wt{A}_\bullet$, both sides of \eqref{RRR-eq} vanish. Thus, the additivity of $\val_{\bZ_p}(-)$ and $\val_{\cO}(-)$ in short exact sequences of complexes and \Cref{Raynaud-re} give \eqref{RRR-eq}. It also follows that the \'{e}tale side of \eqref{RRR-eq} is 
\[
\tst \sum_i (-1)^{i + 1} \ord_p(\#(\Ker f_i)),
\]
which is the negative of the value of the expression 
\[
\tst \sum_i (-1)^i \ord_p\p{\#\p{\f{\Ker \wt{d}_i}{\im \wt{d}_{i - 1}}}} \qq \text{for} \qq (\Ker f_\bullet, \wt{d}_\bullet).
\]
Since this expression is additive in short exact sequences of complexes, the claim about the value of both sides of \eqref{RRR-eq} follows.
\QED

\bremst
\remit \lab{finer-vals}
In \Cref{Raynaud-re-re}, suppose that the $A_i$ are endowed with an action of a commutative ring $\bT$ that is finite free as a $\bZ$-module and that the $d_i$ are $\bT$-equivariant. Then, under the assumption that each $\im d_{i - 1}$ has a $\bT$-stable isogeny complement $B_i \subset A_{i}$, the proof of \Cref{Raynaud-re-re} gives a $\bT$-equivariant conclusion: for a maximal ideal $\fm \subset \bT$ of residue characteristic $p$, we have
\[
\tst\qq  \val_{\bZ_p}\p{H^1_\et((A_\bullet)_{\ov{K}},\, \bZ_p)_\fm} = \val_{\cO}\p{H^1_\dR(\cA_\bullet/\cO)_\fm}.
\]

\remit
The \'{e}tale and the de Rham complexes of \Cref{Raynaud-re-re} are exact after inverting $p$ and perfect, so, as in \cite{KM76}*{Ch.~II, esp.~pp.~47--48}, they have associated Cartier divisors 
\[
\qq\Div(H^1_\et(A_\bullet,\, \bZ_p)) \q \text{on}\q \Spec \bZ_p \qq \text{and} \qq \Div(H^1_\dR(\cA_\bullet/\cO)) \q\text{on} \q \Spec \cO.
\]
Therefore, \cite{KM76}*{Thm.~3~(vi)} reformulates \eqref{RRR-eq} into the following relation between the degrees of those divisors:
\[
\tst e(K/\bQ_p) \cdot \deg_{\bZ_p}(\Div(H^1_\et((A_\bullet)_{\ov{K}}, \bZ_p))) = \deg_{\cO}(\Div(H^1_\dR(\cA_\bullet/\cO))).
\]

\remit
The \'{e}tale side of \eqref{RRR-eq} (or of \eqref{RR-eq}) is invariant under passage to a finite extension of $K$. Thus, even though the N\'{e}ron models may change, the de Rham side is invariant as well. %At least without the $\bT$-action, i.e., in the former case, this may also be deduced from \cite{Cha00}*{6.7}, which, however, also ultimately boils down to \cite{Ray85}*{2.1.1}.
\eremst

\beg \lab{eg-of-A}
By \Cref{Raynaud-re-re}, for a short exact sequence 
\[
0 \ra A \ra B \ra C \ra 0
\]
of abelian varieties over $K$, both sides of \eqref{RRR-eq} vanish. Thanks to the filtrations \eqref{filt}, this vanishing of the de Rham side means that
\[
\tst \length_{\cO}\p{\f{\Lie \cC}{\im(\Lie \cB)}} - \length_{\cO}\p{\f{\Ker(\Lie \cB \ra \Lie \cC)}{\Lie \cA}} = \length_{\cO}\p{\f{\Lie \cA^\vee}{\im(\Lie \cB^\vee)}} - \length_{\cO}\p{\f{\Ker(\Lie \cB^\vee \ra \Lie \cA^\vee)}{\Lie \cC^\vee}}.
\]
It is explained in \cite{LLR04}*{proof of Thm.~2.1} how one associates a smooth finite type $\ov{k}$-group scheme $D$ (resp.,~$D'$) to the morphism $\cB \ra \cC$ (resp.,~$\cB^\vee \ra \cA^\vee$) in such a way that 
\[
D(\ov{k}) \cong \Coker(\cB(\wh{\cO^\sh}) \ra \cC(\wh{\cO^\sh})) \qq \text{ (resp.,} \qq  D'(\ov{k}) \cong \Coker(\cB^\vee(\wh{\cO^\sh}) \ra \cA^\vee(\wh{\cO^\sh}))),
\]
where $\wh{\cO^\sh}$ denotes the completion of the strict Henselization of $\cO$. From this optic, by \cite{LLR04}*{Thm.~2.1~(b)} and \cite{BLR90}*{7.1/6}, \Cref{Raynaud-re-re} proves the equality
\[
\dim D = \dim D'.
\]
\eeg

\Cref{eg-of-A} together with the exactness results \cite{BLR90}*{7.5/4 (ii)} and \cite{AU96}*{Thm.~A.1} of Raynaud leads us to the following corollary (compare with \cite{AU96}*{Cor.~A.2}).

\bcor \lab{dual-exact}
Suppose that $e(K/\bQ_p) \le p - 1$ and let 
\[
0 \ra A \ra B \ra C \ra 0
\]
be a short exact sequence of abelian varieties over $K$ such that $B$ has semiabelian reduction. Then the sequences
\[
0 \ra \Lie \cA \ra \Lie \cB \ra \Lie \cC  \qq \text{and} \qq
0 \ra \Lie \cC^\vee \ra \Lie \cB^\vee \ra \Lie \cA^\vee 
\]
are left exact and both $\f{\Lie \cC}{\im(\Lie \cB)}$ and $\f{\Lie \cA^\vee}{\im(\Lie \cB^\vee)}$ are of the form $(\cO/p\cO)^r$ with the same $r \in \bZ_{\ge 0}$. \QED
\ecor

\brem \lab{finer-dual-exact}
In \Cref{dual-exact}, suppose that $A$, $B$, and $C$ are endowed with an action of a commutative ring $\bT$ that is finite free as a $\bZ$-module, that the sequence is $\bT$-equivariant, and that there is a $\bT$-stable abelian subvariety $C' \subset B$ that maps isogenously to $C$. Then Remark \ref{finer-vals} leads to a further $\bT$-equivariant conclusion: for a maximal ideal $\fm \subset \bT$ of residue characteristic~$p$,
\[
\Coker\p{(\Lie \cB)_\fm \ra (\Lie \cC)_\fm} \simeq \Coker\p{(\Lie \cB^\vee)_\fm \ra (\Lie \cA^\vee)_\fm}.
\]
\erem

\Cref{dual-exact} gives the following consequence of the semistable case of the Manin conjecture. For odd $p$, this consequence also follows from exactness properties of N\'{e}ron models \cite{BLR90}*{7.5/4}.

\bcor \lab{cl-im}
For a new elliptic optimal quotient $\pi \colon J_0(n) \surjects E$ and a prime $p$, if $p^2 \nmid n$, then $\pi$ \up{resp.,~$\pi^\vee$} induces a smooth morphism \up{resp.,~a closed immersion} on N\'{e}ron models over $\bZ_p$ and the sequence
\be \lab{CI-seq}
0 \ra H^1_\dR(\cE/\bZ_p) \ra H^1_\dR(\cJ/\bZ_p) \ra H^1_\dR(\cK/\bZ_p) \ra 0
\ee
is short exact, where $\cE$, $\cJ$, and $\cK$ denote the N\'{e}ron models over $\bZ_p$ of $E$, $J_0(n)$, and $\Ker \pi$.
\ecor

\bpf
\Cref{manin-H-sst} gives the claim about $\pi$. It then follows from \Cref{dual-exact} that the map $J_0(n)^\vee \ra (\Ker \pi)^\vee$ also induces a smooth morphism on N\'{e}ron models over $\bZ_p$. Thus, by \cite{BLR90}*{7.1/6}, the map $\pi^\vee$ induces a closed immersion. Then the Lie algebra complexes that (via \eqref{filt}) comprise the graded pieces of the Hodge filtration of \eqref{CI-seq}  are short exact, so \eqref{CI-seq} must be,~too.
\epf

The proof of the following result illustrates \Cref{Raynaud-re-re} beyond isogenies and short exact sequences.

\bprop \lab{BMS-eg}
Let $A$ be an abelian variety over $K$ endowed with an action of a commutative ring $\bT$ that is finite free as a $\bZ$-module. For every ideal $\fn \subset \bT$ such that $n \in \fn$ for some $n \in \bZ_{\ge 1}$,
\be \lab{spec-ineq}
\tst \val_{\bZ_p}\p{\f{H^1_\et(A_{\ov{K}},\, \bZ_p)}{\fn \cdot H^1_\et(A_{\ov{K}},\, \bZ_p)}} \le \val_{\cO}\p{\f{H^1_\dR(\cA/\cO)}{\fn \cdot H^1_\dR(\cA/\cO)}}.
\ee
\eprop

The following heuristic suggests \eqref{spec-ineq}: one may hope that a suitable formalism of integral $p$-adic Hodge theory would realize the quotient on the right side of \eqref{spec-ineq} as a specialization of the one on the left side, and (normalized) length cannot decrease under specialization. In the good reduction case, one can indeed prove \Cref{BMS-eg} in this way by using the results of \cite{BMS16}.

\bpf
We choose generators $n_1, \ldots, n_m \in \bT$ of $\fn$ and consider the complex of abelian varieties
\[
\tst A \xra{(n_1,\, \ldots,\, n_m)} \prod_{i = 1}^m A \xra{q} Q, \qq \text{where $Q$ is defined to be the cokernel of the first map.}
\]
This complex is exact up to isogeny because the kernel of the first map is killed by $n$. Therefore, \Cref{Raynaud-re-re} applies and gives the following equality, which implies \eqref{spec-ineq}:
\[
\tst \val_{\bZ_p}\p{\f{H^1_\et(A_{\ov{K}}, \bZ_p)}{\fn \cdot H^1_\et(A_{\ov{K}}, \bZ_p)}} = \val_{\cO}\p{\f{H^1_\dR(\cA/\cO)}{\fn \cdot H^1_\dR(\cA/\cO)}} - \val_{\cO}\p{\f{\Ker\p{(n_1,\, \ldots,\, n_m)^*\colon \prod_{i = 1}^m H^1_\dR(\cA/\cO) \ra H^1_\dR(\cA/\cO)}}{q^*(H^1_\dR(\cQ/\cO))}}. \qedhere
\]
% To clearly see that the formation of $H^1_\dR(-/\cO)$ commutes with products, use functoriality to construct a comparison pullback map towards the product of $H^1_\dR$'s and then check that it is an isomorphism by using the filtration to reduce to the compatibility with products of the formation of Lie algebras.
\epf

We end \S\ref{et-dR} by applying \Cref{BMS-eg} in the context of modular curves to exhibit a relation \eqref{Gor-ineq} between the multiplicity of the mod $\fm$ torsion of a modular Jacobian and the Gorenstein defect of the $\fm$-adic completion of the Hecke algebra. In a variety of settings, relations of this sort that are sharper than \eqref{Gor-ineq} are known---see \cite{KW08}*{\S\S1--2} for an overview. Nevertheless, \Cref{Gor}, especially its case \ref{G-ii}, seems to cover some situations that are not addressed in the literature. 

\bcor \lab{Gor}
Fix an $n \in \bZ_{\ge 1}$, let $\bT \subset \End_\bQ(J_0(n))$ be the Hecke algebra defined in \uS\uref{Hecke-def}, let $\fm \subset \bT$ be a maximal ideal of residue characteristic $p$, and let $\cJ$ be the N\'{e}ron model of $J_0(n)$ over $\bZ$. Suppose that $(\Lie \cJ)_\fm$ is free of rank $1$ as a $\bT_\fm$-module, e.g.,~that either of the following holds\ucolon
\benumr
\item \lab{G-i}
$p\nmid n$\uscolon or

\item \lab{G-ii}
$p^2\nmid n$ and $p$ is odd\uscolon or

\item \lab{G-iii}
$\ord_p(n) = 1$ and $U_p \bmod \fm$ lies in $\bF_p^\times \subset \bT/\fm$.
\eenum

Then
\be \lab{Gor-ineq}
\q \dim_{\bT/\fm} (J_0(n)[\fm]) \le \dim_{\bT/\fm} ((\bT/p\bT)[\fm]) + 1
\ee
and $\bT_\fm$ is Gorenstein if and only if $H^1_\et(J_0(n)_{\ov{\bQ}}, \bZ_p)_\fm$ is free \up{of rank $2$} as a $\bT_\fm$-module if and only if $H^1_\dR(\cJ/\bZ_p)_\fm$ is free \up{of rank $2$} as a $\bT_\fm$-module.
\ecor

\bpf
\Cref{H1-free} shows that either of \ref{G-i}--\ref{G-iii} implies the assumed $\bT_\fm$-freeness of $(\Lie \cJ)_\fm$.

Let $w$ denote the Atkin--Lehner involution of $J_0(n)$, so that the induced action of $\bT$ on $J_0(n)^\vee$ is identified with the action of $w\bT w$ on $J_0(n)$ via the isomorphism 
\[
\bT \xra{t \mapsto wtw} w\bT w \qq\text{and the inverse} \qq J_0(n)^\vee \isomto J_0(n)
\]
of the canonical principal polarization (see \cite{MW84}*{Ch.~II, \S5.6 (c)}). The automorphism 
\[
\Lie w \colon \Lie \cJ \isomto \Lie \cJ
\]
intertwines the actions of $\bT$ and of $w\bT w$, so the freeness of $(\Lie \cJ)_\fm$ implies that of $(\Lie \cJ^\vee)_\fm$.\footnote{To stress the Hecke functoriality we let $\cJ^\vee$ denote the N\'{e}ron model of $J_0(n)^\vee$ over $\bZ$, even though $\cJ^\vee \cong \cJ$.} The filtration \eqref{filt} then gives a (necessarily split) extension
\be \lab{dR-T-seq}
 0 \ra \Hom_{\bZ_p}(\bT_\fm, \bZ_p) \ra H^1_\dR(\cJ/\bZ_p)_\fm \ra \bT_\fm \ra 0
 \ee
of $\bT_\fm$-modules, which proves that
\be \lab{dR-mod-m}
\tst \dim_{\bT/\fm} \p{\f{H^1_\dR(\cJ/\bZ_p)}{\fm \cdot H^1_\dR(\cJ/\bZ_p)}} = \dim_{\bT/\fm} ((\bT/p\bT)[\fm]) + 1.
\ee
The modified Weil pairing $(x, y) \mapsto \In{x, wy}_{\mathrm{Weil}}$ shows that $H^1_\et(J_0(n)_{\ov{\bQ}}, \bZ_p)$ is $\bT_{\bZ_p}$-equivariantly isomorphic to its own $\bZ_p$-linear dual (see \cite{DDT97}*{Lem.~1.38}), so 
\be \lab{et-mod-m}
\tst \dim_{\bT/\fm}\p{ \f{H^1_\et(J_0(n)_{\ov{\bQ}},\, \bZ_p)}{\fm\cdot H^1_\et(J_0(n)_{\ov{\bQ}},\, \bZ_p)}} = \dim_{\bT/\fm}(J_0(n)[\fm]).
\ee
The combination of \eqref{dR-mod-m}, \eqref{et-mod-m}, and \Cref{BMS-eg} implies \eqref{Gor-ineq}. If $H^1_\et(J_0(n)_{\ov{\bQ}}, \bZ_p)_\fm$ is free as a $\bT_\fm$-module (necessarily of rank $2$, see \cite{DDT97}*{Lem.~1.39}), then 
\[
\bT_\fm^{\oplus 2} \simeq (\Hom_{\bZ_p}(\bT_\fm, \bZ_p))^{\oplus 2} \q \text{as $\bT_\fm$-modules,}
\] 
so $\bT_\fm$ is Gorenstein. Conversely, if $\bT_\fm$ is Gorenstein, then, due to \eqref{dR-T-seq}, the $\bT_\fm$-module $H^1_\dR(\cJ/\bZ_p)_\fm$ is free of rank $2$, and hence, by \Cref{BMS-eg} and the Nakayama lemma (compare with the paragraph that follows \eqref{s-Nak-arg}), so is $H^1_\et(J_0(n)_{\ov{\bQ}}, \bZ_p)_\fm$.
\epf

\brems
\remi
We are not aware of examples in which the inequality \eqref{Gor-ineq} is sharp.

\remi
\Cref{Gor} leads to examples of non-Gorenstein $\bT_\fm$. For instance, if an odd $n$ is squarefree and $\dim_{\bT/\fm} (J_0(n)[\fm]) > 2$ for some $\fm$ (as happens for $n = 19 \cdot 41$ and an Eisenstein $\fm$ of residue characteristic $5$, see \cite{Yoo16}*{Ex.~4.7}), then, by \Cref{Gor}, the ring $\bT_\fm$ is not Gorenstein.
\erems

%%%%%%%%%%%%%%%%%%%%%%%%%%%%%%%%%%%%%%%

\section{The \'{e}tale and the de Rham congruences} \lab{secondi}

The goal of this section is to use the results of \S\ref{et-dR}, especially \Cref{Raynaud-re}, to derive an exactness result for N\'{e}ron models in \Cref{Ner-exact} in the setting of an abelian variety equipped with a ``Hecke action.'' This exactness result and its applicability criterion given by \Cref{cong-crit} will be useful in \S\ref{contorni} for proving new cases of the generalization of the Manin conjecture to higher dimensional newform quotients. To capture the relevant axiomatics, we begin with an abstract local setup.

\bpp[An abelian variety equipped with rational idempotents] \lab{sec-setup}
Throughout \S\ref{secondi}, we fix
\begin{itemize}
\item
a mixed characteristic $(0, p)$ complete discretely valued field $K$ whose residue field $k$ is perfect;

\item
an abelian variety $A$ over $K$ and its N\'{e}ron model $\cA$ over the ring of integers $\cO$ of $K$;

\item
idempotents $e_1, e_2 \in \End_K(A) \tensor_\bZ \bQ$ that satisfy $e_1 + e_2 = 1$;

\item
a commutative subring $\bT \subset \End_K(A)$ whose elements commute with $e_1$ and $e_2$.

\end{itemize}

As in \S\ref{et-dR}, we let calligraphic letters indicate N\'{e}ron models over $\cO$ (see \S\ref{et-dR-not}).

We are primarily interested in the case when $e_1, e_2 \not\in \End_K(A)$---then various ``rational objects,'' e.g.,~$p$-adic \'{e}tale and de Rham cohomologies, attached to $A$ decompose into summands cut out by the $e_i$, but their integral counterparts typically do not decompose, which produces interesting congruences. As we will see in \S\ref{contorni}, Jacobians of modular curves provide a rich supply of examples of the situation above. Similarly to \Cref{Raynaud-re}, remembering $\bT$ leads to finer statements than with $\bT = \bZ$, and this will be important in applications to newform quotients of modular Jacobians.
\epp

\bpp[The \'{e}tale congruences] \lab{et-cong-numb}
 The \emph{\'{e}tale congruence module} is the finite quotient
\[
\f{H^1_\et(A_{\ov{K}}, \bZ_p)}{H^1_\et(A_{\ov{K}}, \bZ_p)[e_1] + H^1_\et(A_{\ov{K}}, \bZ_p)[e_2]}
\]
(see \eqref{abuse} for the $(-)[e_i]$ notation), and the \emph{\'{e}tale congruence number} is its order. The ring $\bT_{\bZ_p}$ acts, and the \'{e}tale congruence module decomposes compatibly with the decomposition $\bT_{\bZ_p} \cong \prod \bT_\fm$, where $\fm \subset \bT$ ranges over the maximal ideals of residue characteristic $p$.
\epp

\bpp[The de Rham congruences] \lab{dR-cong-numb}
The \emph{de Rham congruence module} is the finite quotient
\[
\f{H^1_\dR(\cA/\cO)}{H^1_\dR(\cA/\cO)[e_1] + H^1_\dR(\cA/\cO)[e_2]}
\]
that we temporarily denote by $M$ (see \S\ref{MM-review} for a review of the lattice $H^1_\dR(\cA/\cO) \subset H^1_\dR(A/K)$), and the \emph{de Rham congruence number} is $p^{\val_\cO(M)}$ with $\val_\cO(\cdot)$ of \S\ref{norm-length}. Like its \'{e}tale counterpart, $M$ decomposes into $\fm$-primary pieces compatibly with $\bT_{\bZ_p} \cong \prod \bT_\fm$.
\epp

We will see in \Cref{Ner-exact} that relations between the \'{e}tale and the de Rham congruence numbers are intricately linked to exactness properties of N\'{e}ron models and in \Cref{cong-differ} that the two numbers differ in general. We begin with a key comparison result.

\bthm \lab{one-ineq}
In the setting of \uS\uref{sec-setup}, for every maximal ideal $\fm \subset \bT$ of residue characteristic $p$,
\[
\tst \val_{\bZ_p}\p{\p{\f{H^1_\et(A_{\ov{K}},\, \bZ_p)}{H^1_\et(A_{\ov{K}},\, \bZ_p)[e_1] + H^1_\et(A_{\ov{K}},\, \bZ_p)[e_2]}
} \tensor_{\bT_{\bZ_p}} \bT_\fm} \ge \val_\cO\p{\p{\f{H^1_\dR(\cA/\cO)}{H^1_\dR(\cA/\cO)[e_1] + H^1_\dR(\cA/\cO)[e_2]}} \tensor_{\bT_{\bZ_p}} \bT_\fm}.
\]
\ethm

The following notation will be useful for the proof of \Cref{one-ineq} and for the subsequent discussion.

\bpp[The abelian varieties $A_i$ and $Q_i$] \lab{Ai-Qi} For an $i \in \{1, 2\}$,
\begin{itemize}
\item
we let $A_i \subset A$ be the image of $A$ under any $\bZ_{> 0}$-multiple of $e_i$ that lies in $\End_K(A)$;

\item
we set $Q_i \ce A/A_i$.
\end{itemize}

Effectively, $A_i$ is the abelian subvariety of $A$ ``cut out by $e_i$'' and is $\bT$-stable, so $A_i$ and $Q_i$ inherit a $\bT$-action. The inclusion and quotient homomorphisms $j_i\colon A_i \hra A$ and  $q_i \colon A \surjects Q_i$ are $\bT$-equivariant.
\epp

\brem \lab{all-Ai}
Any abelian subvariety $B \subset A$ is cut out (in the sense of \S\ref{Ai-Qi}) by some idempotent $e \in \End_K(A) \tensor_\bZ \bQ$. Indeed, this property is isogeny invariant---if $f \colon A \ra A'$ is an isogeny and $e$ cuts out $B$, then $\f{1}{\deg f}\cdot f \circ e \circ f'$ cuts out $f(B)$, where $f' \colon A' \ra A$ is the isogeny such that $f' \circ f = \deg f$---and it clearly holds for the abelian variety $B \times_K A/B$ that is isogenous to $A$.
\erem

\bpf[Proof of Theorem~\uref{one-ineq}]
By construction,
\[
q_i^*(H^1_\et((Q_i)_{\ov{K}}, \bZ_p)) = H^1_\et(A_{\ov{K}}, \bZ_p)[e_{i}],
\]
so \Cref{Raynaud-re} applied to the isogeny $A \xra{q\, =\, (q_1, q_2)} Q_1 \times_K Q_2$ proves that
\[
\tst \val_{\bZ_p}\p{\p{\f{H^1_\et(A_{\ov{K}},\, \bZ_p)}{H^1_\et(A_{\ov{K}},\, \bZ_p)[e_1] + H^1_\et(A_{\ov{K}},\, \bZ_p)[e_2]}
} \tensor_{\bT_{\bZ_p}} \bT_\fm} = \val_\cO \p{ \p{ \f{H^1_\dR(\cA/\cO)}{q_1^*(H^1_\dR(\cQ_1/\cO)) + q_2^*(H^1_\dR(\cQ_2/\cO))} } \tensor_{\bT_{\bZ_p}} \bT_\fm }.
\]
It remains to note that $\f{H^1_\dR(\cA/\cO)}{H^1_\dR(\cA/\cO)[e_1] + H^1_\dR(\cA/\cO)[e_2]}$ is a quotient of $\f{H^1_\dR(\cA/\cO)}{q_1^*(H^1_\dR(\cQ_1/\cO)) + q_2^*(H^1_\dR(\cQ_2/\cO))}$.
\epf

\brem \lab{et-cong-deg}
It follows from the proof and from the last claim of \Cref{Raynaud-re} that the $\fm$-primary factor of the \'{e}tale congruence number equals $\#((A_1 \cap A_2)[\fm^\infty])$.
\erem

\bcor \lab{Ner-exact}
In Theorem~\uref{one-ineq}, if the equality holds\ucolon
\[
\tst \val_{\bZ_p}\p{\p{\f{H^1_\et(A_{\ov{K}},\, \bZ_p)}{H^1_\et(A_{\ov{K}},\, \bZ_p)[e_1] + H^1_\et(A_{\ov{K}},\, \bZ_p)[e_2]}
} \tensor_{\bT_{\bZ_p}} \bT_\fm} = \val_\cO\p{\p{\f{H^1_\dR(\cA/\cO)}{H^1_\dR(\cA/\cO)[e_1] + H^1_\dR(\cA/\cO)[e_2]}} \tensor_{\bT_{\bZ_p}} \bT_\fm},
\]
then the maps $(\Lie \cA)_\fm \ra (\Lie \cQ_i)_\fm$ are surjective. In particular, if the displayed equality holds for every $\fm$, then the $\cO$-morphisms $\cA \ra \cQ_i$ \up{resp.,~$\cA_i \ra \cA$} are smooth \up{resp.,~closed immersions}.
\ecor

\bpf
By the proof of \Cref{one-ineq}, the displayed equality holds if and only if the sequences
\[
0 \ra H^1_\dR(\cQ_i/\cO)_\fm \ra H^1_\dR(\cA/\cO)_\fm \ra H^1_\dR(\cA_i/\cO)_\fm \qq \text{for $i \in \{1, 2\}$}
\]
are left exact. Due to the filtrations \eqref{filt}, this exactness implies the left exactness of the sequences
\[
0 \ra (\Lie \cQ_i)^*_\fm \ra (\Lie \cA)^*_\fm \ra (\Lie \cA_i)^*_\fm,
\]
and hence also the surjectivity of $(\Lie \cA)_\fm \ra (\Lie \cQ_i)_\fm$. For the last claim, it remains to recall that $\cA \ra \cQ_i$ is smooth if and only if $\Lie \cA \ra \Lie \cQ_i$ is surjective and that the smoothness of $\cA \ra \cQ_i$ implies that $\cA_i \ra \cA$ is a closed immersion (see \cite{BLR90}*{7.1/6}).
\epf

\brem \lab{cong-differ}
\Cref{Ner-exact} implies that the \'{e}tale and the de Rham congruence numbers differ in general. Indeed, in the case when $e(K/\bQ_p) \ge p - 1$, there are examples of inclusions $A_1 \subset A$ of abelian varieties (with good reduction) that fail to induce closed immersions on N\'{e}ron models over $\cO$ (see \cite{BLR90}*{7.5/8}) and, by \Cref{all-Ai}, any such $A_1$ arises from some $e_1$ (with $\bT = \bZ$).
\erem

We end \S\ref{secondi} with a criterion for the \'{e}tale and the de Rham congruence numbers to be equal.

\bprop \lab{cong-crit}
In the setup of \uS\uref{sec-setup}, suppose that $\fm \subset \bT$ be a maximal ideal of residue characteristic $p$
for which there exists an $r \in \bZ_{\ge 0}$ such that
\benumr
\item \lab{CC-b}
the $(\bT_\fm \tensor_{\bZ_p} \cO)$-module $H^1_\dR(\cA/\cO)_\fm$ is free of rank $r$\uscolon and

\item \lab{CC-a}
the $(\bT_\fm \tensor_{\bZ_p} \bQ_p)$-module $H^1_\et(A_{\ov{K}}, \bQ_p) \tensor_{\bT_{\bZ_p}} \bT_\fm$ is free of rank $r$.

\eenum

Then the $\bT_\fm$-module $H^1_\et(A_{\ov{K}}, \bZ_p)_\fm$ is free of rank $r$ and the equality holds in Theorem~\uref{one-ineq}\ucolon
\[
\tst \val_{\bZ_p}\p{\p{\f{H^1_\et(A_{\ov{K}},\, \bZ_p)}{H^1_\et(A_{\ov{K}},\, \bZ_p)[e_1] + H^1_\et(A_{\ov{K}},\, \bZ_p)[e_2]}
} \tensor_{\bT_{\bZ_p}} \bT_\fm} = \val_\cO\p{\p{\f{H^1_\dR(\cA/\cO)}{H^1_\dR(\cA/\cO)[e_1] + H^1_\dR(\cA/\cO)[e_2]}} \tensor_{\bT_{\bZ_p}} \bT_\fm}.
\]
\eprop

\bpf
The assumption on $H^1_\dR(\cA/\cO)_\fm$ gives the equality 
\[
\tst \val_{\cO}\p{\f{H^1_\dR(\cA/\cO)}{\fm \cdot H^1_\dR(\cA/\cO)}} = r \cdot \dim_{\bF_p}(\bT/\fm).
\]
Therefore, by \Cref{BMS-eg}, 
\[
\tst \dim_{\bT/\fm} \p{\f{H^1_\et(A_{\ov{K}},\, \bZ_p)}{\fm \cdot H^1_\et(A_{\ov{K}},\, \bZ_p)}} \le r,
\]
to the effect that, by the Nakayama lemma, $H^1_\et(A_{\ov{K}}, \bZ_p)_\fm$ is generated by $r$ elements as a $\bT_\fm$-module. Due to \ref{CC-a}, these $r$ elements are $\bT_\fm$-independent, so the desired $H^1_\et(A_{\ov{K}}, \bZ_p)_\fm \simeq (\bT_\fm)^{\oplus r}$ follows. Consequently, both sides of the claimed equality are equal to 
\[
\tst r \cdot \val_{\bZ_p}\p{\f{\bT_\fm}{\bT_\fm[e_1] + \bT_\fm[e_2]}}. \qedhere
\]
\epf

\brem \lab{Hecke-free}
Modular Jacobians endowed with their Hecke action tend to satisfy \ref{CC-a} (see \cite{DDT97}*{Lem.~1.38--1.39}). Thus, loosely speaking, for them \Cref{cong-crit} proves that the Hecke-freeness of the integral de Rham cohomology implies the Hecke-freeness of the integral $p$-adic \'{e}tale cohomology.
\erem

%%%%%%%%%%%%%%%%%%%%%%%%%%%%%%%%%%%%%%%%%%%%%%%%%%%%%

%%%%%%%%%%%%%%%%%%%%%%%%%%%%%%%%%%%%%%%

\section{The semistable case of the higher dimensional Manin conjecture} \lab{contorni}

The Manin conjecture has been generalized to newform quotients of arbitrary dimension (see \Cref{manin-genl-conj}), and our goal is to address this generalization. More precisely,  we prove in \Cref{genl-counter} that in the higher dimensional case the conjecture fails even at a prime of good reduction and we prove many of its semistable cases in \Cref{genl-results}. Our techniques also supply general relations between the modular degree and the congruence number (see \Cref{deg-cong} and \Cref{DVR-crit,Gor-also-ok}) and determine the endomorphism rings of suitable newform quotients (see \Cref{new-endo}).

\bpp[A newform of level $\Gamma_0(n)$] \lab{mod-case}
Throughout \S\ref{contorni},
\begin{itemize}
\item
we fix an $n \in \bZ_{\ge 1}$ and let $\cJ$ be the N\'{e}ron model over $\bZ$ of $J_0(n)$;

\item
as in \S\ref{Hecke-def}, we let $\bT \subset \End_{\bQ}(J_0(n))$ be the $\bZ$-subalgebra generated by all the $T_\ell$ and $U_\ell$;

\item
we fix a normalized weight $2$ newform $f$ of level $\Gamma_0(n)$, let $e_f \in \bT_\bQ$ be the idempotent that cuts out the factor of $\bT_\bQ$ determined by $f$, and set $e_{f^\perp} \ce 1 - e_f$;

\item
we set $\cO_f \ce \bT/\bT[e_f]$ (see \eqref{abuse}), so that $\cO_f$ is an order in a totally real number field;

\item
we let $\pi_f\colon J_0(n) \surjects A_f$ be the optimal newform quotient determined by $f$, so that 
\[
\q \cO_f \hra \End_\bQ(A_f);
\]

\item
we let $\cA_f$, $\cA_f^\vee$, $\cK$, and $\cK^\vee$ be the N\'{e}ron models over $\bZ$ of $A_f$, $A_f^\vee$, $\Ker \pi_f$, and $(\Ker \pi_f)^\vee$;

\item 
we let $S_2(\Gamma_0(n), \bZ)$ be the module of those weight $2$ cusp forms of level $\Gamma_0(n)$ whose $q$-expansion at the cusp ``$\infty$'' lies in $\bZ\llb q \rrb$, and for a commutative ring $R$ we set 
\[
\q S_2(\Gamma_0(n), R) \ce S_2(\Gamma_0(n), \bZ) \tensor_\bZ R.
\]
\end{itemize}
\epp

The Manin conjecture for $A_f$ is the following generalization of \Cref{manin-conj}.

\bconj[\cite{Joy05}*{Conj.~2} or \cite{ARS06}*{Conj.~3.12}] \lab{manin-genl-conj}
In the setting of \uS\uref{mod-case} \up{see also \eqref{abuse}},
\[
\pi_f^*(H^0(\cA_f, \Omega^1)) = S_2(\Gamma_0(n), \bZ)[e_{f^\perp}] \qq \text{inside} \qq H^0(J_0(n), \Omega^1) \cong S_2(\Gamma_0(n), \bQ).
\]
\econj

\brem \lab{genl-conj-exp}
As in \cite{Edi91}*{proof of Prop.~2}, since ``$\infty$'' extends to a $\bZ$-point of $X_0(n)^\sm$, the N\'{e}ron property gives the inclusion 
\[
\tst \pi_f^*(H^0(\cA_f, \Omega^1)) \subset S_2(\Gamma_0(n), \bZ)[e_{f^\perp}],
\]
so \Cref{manin-genl-conj} amounts to the vanishing of the finite quotient 
\[
\tst \f{S_2(\Gamma_0(n),\, \bZ)[e_{f^\perp}]}{\pi_f^*(H^0(\cA_f,\, \Omega^1))}.
\]
Due to the following standard lemma and the exactness property \cite{BLR90}*{7.5/4 (ii) and its proof} of N\'{e}ron models, the $p$-primary part of this quotient vanishes for every odd prime $p$ with $p^2 \nmid n$ (compare with \cite{ARS06}*{Cor.~3.7}). We will see in \Cref{genl-counter} that, in contrast, the $2$-primary part need not vanish even when $2 \nmid n$.
\erem

\blem \lab{AL-lemma}
In the setting of \uS\uref{mod-case}, if $p$ is a prime with $p \nmid n$, then
\[
S_2(\Gamma_0(n), \bZ_p) = H^0(\cJ_{\bZ_p}, \Omega^1) \qq \text{inside} \qq H^0(J_0(n)_{\bQ_p}, \Omega^1) \cong S_2(\Gamma_0(n), \bQ_p);
\]
 if $p$ is a prime with $p^2 \nmid n$, then 
\[
S_2(\Gamma_0(n), \bZ_p)[e_{f^\perp}] = H^0(\cJ_{\bZ_p}, \Omega^1)[e_{f^\perp}] \qq \text{inside} \qq H^0(J_0(n)_{\bQ_p}, \Omega^1) \cong S_2(\Gamma_0(n), \bQ_p).
\]
\elem

\bpf
If $p \nmid n$, then, due to \cite{Edi06}*{2.5} and \eqref{transfer-T}, 
\[
S_2(\Gamma_0(n), \bZ_p) = H^0(\cJ_{\bZ_p}, \Omega^1).
\] 
By loc.~cit.,~if $p \mid n$ but $p^2 \nmid n$, then 
\[
S_2(\Gamma_0(n), \bZ_p) = H^0(U^\infty, \Omega^1),
\]
where $U^\infty \subset X_0(n)_{\bZ_p}$ is the open complement of the irreducible component of $X_0(n)_{\bF_p}$ that does not contain the reduction of the cusp ``$\infty$.'' Thus, since the $p$-Atkin--Lehner involution $w_p$ interchanges the two irreducible components of $X_0(n)_{\bF_p}$ and acts as $\pm 1$ on $S_2(\Gamma_0(n), \bZ_p)[e_{f^{\perp}}]$, we get 
\[
S_2(\Gamma_0(n), \bZ_p)[e_{f^\perp}] = H^0(X_0(n)_{\bZ_p}, \Omega)[e_{f^\perp}]
\]
(see \cite{Ces16g}*{proof of Lem.~2.7}). The claim then follows from another application of \eqref{transfer-T}.
\epf

Similarly to the elliptic curve case discussed in \S\ref{ell-curve}, the strategy of our analysis of \Cref{manin-genl-conj} is to relate it to a comparison of the congruence number and the modular degree of $f$. We use the following standard lemma to introduce these numbers in \Cref{cong-deg-def}.

\blem \lab{degf-in-N} \hfill
\benum
\item \lab{DIN-a}
The composition $ A_f^\vee \xra{\pi_f^\vee} J_0(n)^\vee \cong J_0(n)$ is $\bT$-equivariant. 

\item \lab{DIN-b}
The finite $\bQ$-group scheme $\pi_f^\vee(A_f^\vee) \cap \Ker \pi_f$ carries a perfect alternating bilinear pairing for which the action of $\bT$ is self adjoint.
\eenum
\elem

\bpf \hfill
\benum
\item
The Atkin--Lehner involution $w$ of $J_0(n)$ acts as $\pm 1$ on $A_f$. Moreover, under the canonical principal polarization $\theta \colon J_0(n) \isomto J_0(n)^\vee$, the action of a $t \in \bT$ corresponds to that of $wt^\vee w$ (see \cite{MW84}*{Ch.~II, \S5.6 (c)}), so the claimed $\bT$-equivariance follows.

\item
Let $\gL \ce \pi_f \circ \theta\i \circ \pi_f^\vee$ be the pullback of $\theta\i$ to a polarization of $A_f^\vee$, so that $A^\vee_f[\gL]$ caries a perfect alternating bilinear Weil pairing (see \cite{Pol03}*{\S10.4, esp.~Prop.~10.3}). Since $\gL$ is $\bT$-equivariant, \cite{Oda69}*{Cor.~1.3~(ii)} ensures that the action of any $t \in \bT$ is self adjoint with respect to this pairing. It remains to observe that 
\[
\q A^\vee_f[\gL] \cong \pi_f^\vee(A_f^\vee) \cap \Ker \pi_f. \qedhere
\]
\eenum
\epf

\bd[Compare with \cite{ARS12}*{\S3}] \lab{cong-deg-def}
The \emph{congruence number} $\mathrm{cong}_f$ of $f$ and the \emph{modular degree} $\deg_f$ of $f$ are the positive integers (see \eqref{abuse} and \Cref{degf-in-N}~\ref{DIN-b})
\[
\tst \mathrm{cong}_{f} \ce \#\p{\f{\bT}{\bT[e_f] + \bT[e_{f^\perp}]}} \qq \text{and} \qq \deg_f \ce \p{ \#\p{\pi_f^\vee(A_f^\vee) \cap \Ker \pi_f}}^{\f{1}{2}};
\]
for a maximal ideal $\fm \subset \bT$, the \emph{$\fm$-primary factors} $\mathrm{cong}_{f,\, \fm}$ and $\deg_{f,\, \fm}$ are the positive integers
\[
\tst \mathrm{cong}_{f,\, \fm} \ce \#\p{\f{\bT_\fm}{\bT_\fm[e_f] + \bT_\fm[e_{f^\perp}]}} \qq \text{and} \qq \deg_{f,\, \fm} \ce \p{\#\p{\p{\pi_f^\vee(A_f^\vee) \cap \Ker \pi_f}[\fm^\infty]}}^{\f{1}{2}}.
\]
\ed

\brem \lab{mod-deg-et-cong}
One may view $\deg_{f,\, \fm}^2$ as (a factor of) an \'{e}tale congruence number: the abelian subvariety of $J_0(n)$ cut out by $e_f$ (resp.,~by $e_{f^\perp}$) as in \S\ref{Ai-Qi} is identified with $\pi_f^\vee(A_f^\vee)$ (resp.,~$\Ker \pi_f$), so
\be \lab{et-rewrite}
\tst \q\deg_{f,\, \fm}^2 \overset{\ref{et-cong-deg}}{=}  \#\p{\f{H^1_\et(J_0(n)_{\ov{\bQ}},\, \bZ_p)_\fm}{H^1_\et(J_0(n)_{\ov{\bQ}},\, \bZ_p)_\fm[e_f] + H^1_\et(J_0(n)_{\ov{\bQ}},\, \bZ_p)_\fm[e_{f^\perp}]}} \q \text{with} \q p \ce \Char(\bT/\fm).
\ee
In particular,
\[
\tst\q\deg_f^2 = \displaystyle\prod_{\text{primes }p} \tst\# \p{\f{H^1_\et(J_0(n)_{\ov{\bQ}},\, \bZ_p)}{H^1_\et(J_0(n)_{\ov{\bQ}},\, \bZ_p)[e_f] + H^1_\et(J_0(n)_{\ov{\bQ}},\, \bZ_p)[e_{f^\perp}]}}.
\]
\erem

Agashe, Ribet, and Stein proved in \cite{ARS12}*{Thm.~3.6} that the exponent of $\pi_f^\vee(A_f^\vee) \cap \Ker \pi_f$ divides the exponent of $\f{\bT}{\bT[e_f] + \bT[e_{f^\perp}]}$ and that both exponents have the same $p$-adic valuation for every prime $p$ with $p^2 \nmid n$. As a key step towards the semistable case of \Cref{manin-genl-conj}, we wish to complement these results with relations between $\mathrm{cong}_f$ and $\deg_f$ themselves. 

\bprop \lab{cong-deg-div}
If $\fm \subset \cO_f$ is a maximal ideal of residue characteristic $p$ with $p^2 \nmid n$, then
\be \lab{CDD-eq}
\deg_{f,\, \fm} = \mathrm{cong}_{f,\, \fm} \cdot \# \Coker\p{(\Lie \cJ)_\fm \ra (\Lie \cA_f)_\fm}.
\ee
\eprop

\bpf
\Cref{dual-exact} gives the left exact sequences
\[
0 \ra (\Lie \cK)_{\bZ_p} \ra (\Lie \cJ)_{\bZ_p} \ra (\Lie \cA_f)_{\bZ_p} \q \text{and} \q 0 \ra (\Lie \cA^\vee_f)_{\bZ_p} \ra (\Lie \cJ)_{\bZ_p} \ra (\Lie \cK^\vee)_{\bZ_p},
\]
which give the following identifications: 
\[
\,\,\,\, (\Lie \cK)_{\bZ_p} \cong (\Lie \cJ)_{\bZ_p}[e_{f}] \q  \text{and} \q (\Lie \cA^\vee_f)_{\bZ_p} \cong (\Lie \cJ)_{\bZ_p}[e_{f^\perp}].
\]
Therefore, \Cref{Raynaud-re} and \eqref{RR-diag} applied to the $\bT$-equivariant (see \Cref{degf-in-N}~\ref{DIN-a}) isogeny $\Ker \pi_f \times A_f^\vee \ra J_0(n)$ give the equality
\be \lab{RR-application}
\tst \deg_{f,\, \fm}^2 = \#\p{\p{\f{\Lie \cJ}{(\Lie \cJ)[e_f] + (\Lie \cJ)[e_{f^\perp}]}
} \tensor_\bT \bT_\fm} \cdot \#\p{\p{\f{\Lie (\cK^\vee)\, \oplus\, \Lie(\cA_f)}{\Lie \cJ}
} \tensor_\bT \bT_\fm}.
\ee
By \Cref{H1-free}~\ref{H1f-i}--\ref{H1f-iii}, we have $(\Lie \cJ)_\fm \simeq \bT_\fm$, so the first factor on the right side of \eqref{RR-application} equals $\mathrm{cong}_{f,\, \fm}$. In addition, by \Cref{finer-dual-exact},
\[
\# \Coker\p{(\Lie \cJ)_\fm \ra (\Lie \cA_f)_\fm} = \# \Coker\p{(\Lie \cJ)_\fm \ra (\Lie \cK^\vee)_\fm},
\]
so, since the kernels of these maps are $(\Lie \cJ)_\fm[e_f]$ and $(\Lie \cJ)_\fm[e_{f^\perp}]$, respectively, the second factor on the right side of \eqref{RR-application} equals
\[
\tst \p{\# \Coker\p{(\Lie \cJ)_\fm \ra (\Lie \cA_f)_\fm}}^2 \cdot \# \Coker\p{ (\Lie \cJ)_\fm \hra \f{(\Lie \cJ)_\fm}{(\Lie \cJ)_\fm[e_{f}]} \oplus  \f{(\Lie \cJ)_\fm}{(\Lie \cJ)_\fm[e_{f^\perp}]}}. 
\]
It remains to observe the short exact sequence
\[
\tst 0 \ra (\Lie \cJ)_\fm \ra \f{(\Lie \cJ)_\fm}{(\Lie \cJ)_\fm[e_{f}]} \oplus  \f{(\Lie \cJ)_\fm}{(\Lie \cJ)_\fm[e_{f^\perp}]} \xra{(x, y)\, \mapsto\, x - y} \f{(\Lie \cJ)_\fm}{(\Lie \cJ)_\fm[e_{f}] + (\Lie \cJ)_\fm[e_{f^\perp}]} \ra 0
\]
and recall that $(\Lie \cJ)_\fm \simeq \bT_\fm$.
\epf

\bcor  \lab{deg-cong}
If $p$ is a prime with $p^2 \nmid n$, then
\be \lab{deg-cong-eq}
\ord_p(\deg_f) = \ord_p(\mathrm{cong}_f) + \ord_p\p{\# \Coker\p{\Lie \cJ \ra \Lie \cA_f}};
\ee
in particular, if, in addition, $p$ is odd, then $\ord_p(\deg_f) = \ord_p(\mathrm{cong}_f)$.
\ecor

\bpf
The product of the equalities \eqref{CDD-eq} over all $\fm$ of residue characteristic $p$ gives \eqref{deg-cong-eq}. In the case when $p$ is odd, \cite{BLR90}*{7.5/4 (ii) and its proof} ensure the smoothness of $\cJ_{\bZ_p} \ra (\cA_f)_{\bZ_p}$, so the second summand of the right side of \eqref{deg-cong-eq} vanishes.
\epf

With \Cref{deg-cong} in hand, we are ready to present a counterexample to \Cref{manin-genl-conj}.

\bthm \lab{genl-counter}
Conjecture~\uref{manin-genl-conj} fails for a $24$-dimensional optimal newform quotient of $J_0(431)$ and also for a $91$-dimensional optimal newform quotient of $J_0(2089)$ \up{both $431$ and $2089$ are primes}.
\ethm

\bpf
By \cite{ARS06}*{Rem.~3.7}, there is a weight $2$ newform $f$ of level $\Gamma_0(431)$ with $\dim A_f = 24$ and 
\[
\deg_f = 2^{11} \cdot 6947 \qq  \text{and} \qq \mathrm{cong}_f = 2^{10} \cdot 6947.
\]
The following Sage code \cite{Sage} confirms these values $\deg_f$ and $\mathrm{cong}_f$.
\verbinput{verb.txt}
Alternatively, the following Magma code \cite{Magma} computes $\deg_f^2$ with a faster runtime.
\verbinput{verb2.txt}
These means also give us a weight $2$ newform $f'$ of level $\Gamma_0(2089)$ with $\dim A_{f'} = 91$ and
\[
\qqq \deg_{f'} = 2^{80} \cdot 3 \cdot 5 \cdot 11 \cdot 19 \cdot 73 \cdot 139 \qq  \text{and} \qq \mathrm{cong}_{f'} = 2^{79} \cdot 3 \cdot 5 \cdot 11 \cdot 19 \cdot 73 \cdot 139.
\]
By \Cref{deg-cong}, the map 
\[
(\Lie \cJ)_{\bZ_2} \ra (\Lie \cA_f)_{\bZ_2}
\]
is not surjective, to the effect that $\f{H^0(\cJ_{\bZ_2},\, \Omega^1)}{\pi_f^*(H^0((\cA_f)_{\bZ_2},\, \Omega^1))}$ has nontrivial $2$-torsion. Since, by \Cref{AL-lemma}, 
\[
H^0(\cJ_{\bZ_2}, \Omega^1) = S_2(\Gamma_0(431), \bZ_2) \qq \text{and} \qq \tst \f{S_2(\Gamma_0(431),\, \bZ_2)}{S_2(\Gamma_0(431),\, \bZ_2)[e_{f^\perp}]} \q \text{is torsion free,}
\]
we get that 
\[
\pi_f^*(H^0(\cA_f, \Omega^1)) \neq S_2(\Gamma_0(431),\, \bZ)[e_{f^\perp}]
\]
(and likewise for $f'$), contrary to \Cref{manin-genl-conj}.
\epf

\brems
\remi \lab{Kil-plug}
\Cref{Gor-also-ok} below and \cite{Kil02} suggest considering the levels $\Gamma_0(431)$, $\Gamma_0(503)$, and $\Gamma_0(2089)$.

\remi \lab{CES-counter}
For the $f$ considered in the proof of \Cref{genl-counter}, let $\wt{\pi}_f \colon J_1(431) \surjects \wt{A}_f$ be the resulting optimal newform quotient of $J_1(431)$. As in the proof of \Cref{down}, there is a commutative diagram
\[
\xymatrix{
J_1(431) \ar@{->>}[r]^-{\wt{\pi}_f} \ar@{->>}[d] & \wt{A}_f \ar@{->>}[d]^-{a} \\
J_0(431) \ar@{->>}[r]^-{\pi_f} & A_f
}
\]
in which $a$ is an isogeny and $\Ker a$ is a quotient of the Cartier dual of the Shimura subgroup $\Sigma(431) \subset J_0(431)$. Moreover, by \cite{LO91}*{Cor.~1 to Thm.~1}, we have $2 \nmid \# \Sigma(431)$, so $\Lie a$ is an isomorphism on Lie algebras of N\'{e}ron models over $\bZ_2$, and hence $\Lie \wt{\pi}_f$ is not surjective on such Lie algebras. Then, as in the proof of \Cref{genl-counter}, the quotient $\f{S_2(\Gamma_1(431),\, \bZ)}{(\wt{\pi}_f)^*(H^0(\wt{\cA}_f,\, \Omega^1))}$ has nontrivial $2$-torsion, where $\wt{\cA}_f$ denotes the N\'{e}ron model  over $\bZ$ of $\wt{A}_f$ and \cite{CES03}*{Lem.~6.1.6} supplies the inclusion 
\[
\qq (\wt{\pi}_f)^*(H^0(\wt{\cA}_f, \Omega^1)) \subset S_2(\Gamma_1(431), \bZ).
\]
This is a counterexample to the analogue \cite{CES03}*{Conj.~6.1.7} of \Cref{manin-genl-conj} for newform quotients of $J_1(n)$.
\erems

In the case when $A_f$ is an elliptic curve, one knows that $\deg_f \mid \mathrm{cong}_f$ (see \Cref{deg-cong-div}). Even though this divisibility fails in the higher dimensional case (see the proof of \Cref{genl-counter}), we wish to generalize it as follows.

\bprop \lab{converse-div}
For a maximal ideal $\fm \subset \cO_f$ of residue characteristic $p$ such that $(\cO_f)_\fm$ is a discrete valuation ring, $\deg_{f,\, \fm} \mid \mathrm{cong}_{f,\, \fm}$.  In particular, if the order $\cO_f$ is maximal, then $\deg_f \mid \mathrm{cong}_f$.
\eprop

\bpf
It will suffice to mildly generalize the proof of \Cref{deg-cong-div}. Namely, 
\[
\tst H^1_\et(J_0(n)_{\ov{\bQ}}, \bQ_p) \tensor_{\bT_{\bZ_p}} \bT_\fm \q \text{is a free $\bT_\fm[\f{1}{p}]$-module of rank $2$}
\]
(see \cite{DDT97}*{Lem.~1.38--1.39}) and $\cO_f = \bT/\bT[e_{f}]$, so, since $(\cO_f)_\fm$ is a discrete valuation ring, we check over $(\cO_f)_{\fm}[\f{1}{p}]$~that
\[
\tst \f{H^1_\et(J_0(n)_{\ov{\bQ}},\, \bZ_p)_\fm}{H^1_\et(J_0(n)_{\ov{\bQ}},\, \bZ_p)_\fm[e_{f}]} \qq \text{is a free $(\cO_f)_\fm$-module of rank $2$.}
\]
Then the further quotient 
\[
\tst \f{H^1_\et(J_0(n)_{\ov{\bQ}},\, \bZ_p)_\fm}{H^1_\et(J_0(n)_{\ov{\bQ}},\, \bZ_p)_\fm[e_{f}] + H^1_\et(J_0(n)_{\ov{\bQ}},\, \bZ_p)_\fm[e_{f^\perp}]} \qq \text{admits a surjection from} \qq \p{\f{\bT_\fm}{\bT_\fm[e_f] + \bT_\fm[e_{f^\perp}]}}^2,
\]
and the desired conclusion follows from \eqref{et-rewrite}.
\epf

\brem
%\remi
%The proof of \Cref{converse-div} continues to work in other contexts. For example, for every subgroup $H$ with $\Gamma_1(n) \subset H \subset \Gamma_0(n)$, it gives an analogous conclusion for newforms of level $H$ (one also includes the diamond operators in the definition of $\bT$ when $H \neq \Gamma_0(n)$).
\Cref{converse-div} implies that for the newform $f$ used in the proof of \Cref{genl-counter}, the order $\cO_f$ is nonmaximal. Indeed, according to \cite{Ste99}*{Table 2}, $\cO_f$ has index $4$ in its normalization.
\erem

As we will see in \Cref{genl-results}, the following consequence of the work above proves the semistable case of \Cref{manin-genl-conj} under the assumption that $\cO_f \tensor_{\bZ} \bZ_2$ is regular. Furthermore, its freeness aspect supplies additional information that seems new in the case of an odd $p$.

\bthm \lab{DVR-crit}
If $\fm \subset \cO_f$ is a maximal ideal of residue characteristic $p$ with $p^2 \nmid n$ such that either $p$ is odd or $(\cO_f)_\fm$ is a discrete valuation ring, then 
\[
\tst (\Lie \cJ)_\fm \ra (\Lie \cA_f)_\fm
\]
is surjective, $\deg_{f,\, \fm} = \mathrm{cong}_{f,\,\fm}$, and $(\Lie \cA_f)_\fm$ is free of rank $1$ as an $(\cO_f)_\fm$-module. 
\ethm

\bpf
In both cases, by \Cref{H1-free}~\ref{H1f-i}--\ref{H1f-iii}, 
\[
(\Lie \cJ)_\fm \simeq \bT_\fm \qq \text{as $\bT_\fm$-modules.}
\]
If $p$ is odd, then the surjectivity of $(\Lie \cJ)_{\fm} \ra (\Lie \cA_f)_\fm$ follows from the smoothness of the map $\cJ_{\bZ_p} \ra (\cA_f)_{\bZ_p}$ supplied by \cite{BLR90}*{7.5/4 (ii) and its proof}. If $(\cO_f)_\fm$ is a discrete valuation ring, then the surjectivity follows by combining \Cref{cong-deg-div,converse-div}. Therefore, in both cases $\deg_{f,\, \fm} \overset{\ref{cong-deg-div}}{=} \mathrm{cong}_{f,\,\fm}$ and $(\Lie \cA_f)_\fm$ is isomorphic to 
\[
\bT_\fm/\bT_\fm[e_{f}] \cong (\cO_f)_\fm
\]
as a $\bT_\fm$-module (so also as an $(\cO_f)_\fm$-module).
\epf

Our next aim is to show in \Cref{Gor-also-ok} that conclusions like those of \Cref{DVR-crit} may also be drawn if the assumption on $(\cO_f)_\fm$ is replaced by the Gorensteinness of $\bT_\fm$. To put the latter condition into context, we now review some cases in which it is known to hold.

\bpp[Multiplicity one] \lab{mult-one}
In the setting of \S\ref{mod-case}, suppose that $p^2\nmid n$, let $\fm \subset \cO_f$ be a maximal ideal of residue characteristic $p$, and let 
\[
\rho_\fm\colon \Gal(\ov{\bQ}/\bQ) \ra \GL_2(\cO_f/\fm)
\]
be the associated semisimple modular mod $\fm$ Galois representation. By \Cref{Gor} and \eqref{et-mod-m}, 
\[
\bT_\fm \qq \text{is Gorenstein if and only if} \qq \dim_{\bT/\fm}(J_0(n)[\fm]) = 2.
\]
Therefore, $\bT_\fm$ is Gorenstein in any of the following cases: \benuma
\item \lab{mult-one-1}
if $\ord_p(n) = 0$ and $\rho_\fm$ is absolutely irreducible, except possibly when, in addition, $p = 2$, the ideal $\fm$ contains $T_2 - 1$, the restriction of $\rho_\fm$ to $\Gal(\ov{\bQ}_2/\bQ_2)$ is unramified, and $\rho_\fm(\Frob_2)$ lies in the center of $\GL_2(\bT/\fm)$, see \cite{Rib90}*{Thm.~5.2~(b)}, \cite{Edi92b}*{Thm.~9.2} (and possibly also \cite{Gro90}*{Thm.~12.10~(1)}), and \cite{RS01}*{Thm.~6.1} in the appendix by Buzzard;

\item
if $\ord_p(n) = 1$ and $\rho_\fm$ is absolutely irreducible and not of level $n/p$, see \cite{MR91}*{Main Thm.};

\item
if $\ord_p(n) = 1$ with $p$ odd, $\rho_\fm$ is irreducible, $U_p \not\in \fm$, and the semisimplification 
\[
\q \tst (\rho_\fm|_{\Gal(\ov{\bQ}_p/\bQ_p)})^{\mathrm{ss}}
\]
is not of the form $\chi \oplus \chi$ for some character $\chi \colon \Gal(\ov{\bQ}_p/\bQ_p) \ra (\bT/\fm)^\times$, see \cite{Wil95}*{Thm.~2.1~(ii)}.
\eenum
\epp

In addition, in the case when $n$ is a prime, $\bT_\fm$ is also Gorenstein when $\fm$ is Eisenstein, see \cite{Maz77}*{16.3}. Also, in contrast to \ref{mult-one-1}, if $p = 2$, then even  when $n$ is a prime $\bT_\fm$ need not be Gorenstein---see \cite{Kil02}; alternatively, the possible failure of Gorensteinness in such situations may be deduced by combining the examples in the proof of \Cref{genl-counter} with the following result.

\bthm \lab{Gor-also-ok}
If $\fm \subset \cO_f$ is a maximal ideal of residue characteristic $p$ with $p^2 \nmid n$ such that $\bT_\fm$ is Gorenstein \up{see \uS\uref{mult-one}}, then 
\[
\tst (\Lie \cJ)_\fm \ra (\Lie \cA_f)_\fm
\]
is surjective, $\deg_{f,\,\fm} = \mathrm{cong}_{f,\,\fm} $, the $(\cO_f)_\fm$-modules $(\Lie \cA_f)_\fm$ and $(\Lie \cA_f^\vee)^*_\fm$ are free of rank $1$, where $(-)^*$ denotes the $\bZ_p$-linear dual, and the $(\cO_f)_\fm$-module $H^1_\dR(\cA_f^\vee/\bZ_p)_\fm$ is free of rank $2$.
\ethm

\bpf
By \Cref{H1-free}~\ref{H1f-i}--\ref{H1f-iii}, 
\[
(\Lie \cJ)_\fm \simeq \bT_\fm \q \text{as $\bT_\fm$-modules,}
\]
so, since $\bT_\fm$ is Gorenstein, also 
\[
(\Lie \cJ)^*_\fm \simeq \bT_\fm \q \text{as $\bT_\fm$-modules}
\]
and, by the proof of \Cref{Gor}, 
\[
(\Lie \cJ^\vee)_\fm \simeq \bT_\fm \q \text{as $\bT_\fm$-modules.\footnote{To stress the Hecke functoriality we let $\cJ^\vee$ denote the N\'{e}ron model of $J_0(n)^\vee$ over $\bZ$, even though $\cJ^\vee \cong \cJ$.}}
\]
The filtration \eqref{filt} then proves that $H^1_\dR(\cJ/\bZ_p)_\fm$ is free of rank $2$ as a $\bT_\fm$-module. Therefore, \Cref{cong-crit} and \Cref{Hecke-free} imply that the $\fm$-primary parts of the \'{e}tale and the de Rham congruence numbers of $J_0(n)_{\bQ_p}$ formed with respect to the idempotents $e_f$ and $e_{f^\perp}$ are equal. Thus, since the abelian subvariety of $J_0(n)$ cut out by $e_f$ (resp.,~by $e_{f^\perp}$) as in \S\ref{Ai-Qi} is identified with $\pi_f^\vee(A_f^\vee)$ (resp.,~with $\Ker \pi_f$), \Cref{Ner-exact} proves the surjectivity of the maps 
\[
(\Lie \cJ)_\fm \ra (\Lie \cA_f)_\fm \qq \text{and} \qq (\Lie \cJ^\vee)_\fm \ra (\Lie \cK^\vee)_\fm.
\]
\Cref{cong-deg-div} then gives $\deg_{f,\,\fm} = \mathrm{cong}_{f,\,\fm}$ and \Cref{dual-exact} gives the short exact sequences
\[
0 \ra (\Lie \cK)_\fm \ra (\Lie \cJ)_\fm \ra (\Lie \cA_f)_\fm \ra 0 \q \text{and} \q 0 \ra (\Lie \cK^\vee)_\fm^* \ra  (\Lie \cJ^\vee)_\fm^* \ra (\Lie \cA_f^\vee)_\fm^* \ra 0,
\]
which show that 
\[
(\Lie \cK)_\fm \cong (\Lie \cJ)_\fm[e_{f}] \qq \text{and}\qq (\Lie \cK^\vee)_\fm^* \cong (\Lie \cJ^\vee)_\fm^*[e_{f}].
\]
Since $\bT_\fm$ is Gorenstein, $(\Lie \cJ^\vee)_\fm^*$ inherits $\bT_\fm$-freeness from $(\Lie \cJ^\vee)_\fm$, and it follows that $(\Lie \cA_f)_\fm$ and $(\Lie \cA_f^\vee)_\fm^*$ are free of rank $1$ as $(\cO_f)_\fm$-modules. The filtration \eqref{filt} then gives the claim about $H^1_\dR(\cA_f^\vee/\bZ_p)_\fm$.
\epf

The following corollary supplies information about endomorphism rings of newform quotients of $J_0(n)$  (we leave the explication of its ``one $\fm$ at a time'' generalization to an interested reader).

\bcor \lab{new-endo}
If $n$ is squarefree and for every maximal ideal $\fm \subset \cO_f$ of residue characteristic $2$ either $(\cO_f)_\fm$ is a discrete valuation ring or $\bT_\fm$ is Gorenstein, then the inclusion $\cO_f \hra \End_\bQ(A_f)$ is an isomorphism.
\ecor

\bpf
\Cref{DVR-crit,Gor-also-ok} imply that the $\cO_f$-module $\Lie \cA_f$ is locally free of rank $1$, to the effect that the inclusion $\cO_f \hra \End_{\cO_f}(\Lie \cA_f)$ is an isomorphism. Since the action of $\End_\bQ(A_f)$ on $\Lie \cA_f$ is faithful and $\cO_f$-linear, the desired conclusion follows.
\epf

We are ready to put the preceding results together to prove suitable semistable cases of \Cref{manin-genl-conj}.

\bthm \lab{genl-results}
In the setting of \uS\uref{mod-case}, suppose that $p$ is a prime with $p^2 \nmid n$ such that either
\benumr
\item
$p$ is odd, or 

\item
for every maximal ideal $\fm \subset \cO_f$ of residue characteristic $p$, either $(\cO_f)_\fm$ is a discrete valuation ring or $\bT_\fm$ is Gorenstein.
\eenum
Then the $p$-primary aspect of Conjecture~\uref{manin-genl-conj} holds for $A_f$ in the sense that 
\[
\tst p \nmid \# \p{\f{S_2(\Gamma_0(n),\, \bZ)[e_{f^\perp}]}{\pi_f^*(H^0(\cA_f,\, \Omega^1))}}
\]
\up{see Remark~\uref{genl-conj-exp}}\uscolon in addition, the map $\pi_f \colon \cJ_{\bZ_p} \ra (\cA_f)_{\bZ_p}$ is smooth, the map $\pi_f^\vee \colon (\cA_f^\vee)_{\bZ_p} \ra (\cJ)_{\bZ_p}$ is a closed immersion, and the sequence
\be \lab{GR-seq}
0 \ra H^1_\dR(\cA_f/\bZ_p) \ra H^1_\dR(\cJ/\bZ_p) \ra H^1_\dR(\cK/\bZ_p) \ra 0
\ee
is short exact.
\ethm

\bpf
By \Cref{DVR-crit,Gor-also-ok}, the map
\[
\tst (\Lie \cJ)_{\bZ_p} \ra (\Lie \cA_f)_{\bZ_p}
\]
is surjective, so $\cJ_{\bZ_p} \ra (\cA_f)_{\bZ_p} $ is smooth and $\f{H^0(\cJ_{\bZ_p}, \,\Omega^1)}{\pi_f^*(H^0((\cA_f)_{\bZ_p},\, \Omega^1))}$ is $p$-torsion free, to the effect that 
\[
\tst \pi_f^*(H^0((\cA_f)_{\bZ_p}, \Omega^1)) = H^0(\cJ_{\bZ_p}, \Omega^1)[e_{f^\perp}].
\]
Since 
\[
H^0(\cJ_{\bZ_p}, \Omega^1)[e_{f^\perp}] \overset{\ref{AL-lemma}}{=} S_2(\Gamma_0(n), \bZ_p)[e_{f^\perp}],
\]
the claim about \Cref{manin-genl-conj} follows. The additional claims then follow from \Cref{dual-exact}, similarly to the proof of \Cref{cl-im}.
\epf

\brem
The short exactness of the sequence \eqref{GR-seq} fails in the counterexamples of \Cref{genl-counter} when $p = 2$. Indeed, due to \eqref{filt}, such exactness would imply the surjectivity of the map
\[
(\Lie \cJ)_{\bZ_2} \ra (\Lie \cK^\vee)_{\bZ_2},
\]
that is, by \Cref{dual-exact}, the surjectivity of the map
\[
(\Lie \cJ)_{\bZ_2} \ra (\Lie \cA_f)_{\bZ_2}.
\]
\erem

In \cite{CES03}*{Conj.~6.1.7}, Conrad, Edixhoven, and Stein proposed an analogue of \Cref{manin-genl-conj} for newform quotients of $J_1(n)$. We now deduce several special cases of their conjecture from \Cref{genl-results} (Remark \ref{CES-counter} exhibited a counterexample to the general case).

\bcor \lab{CES-positive}
Let $p$ be a prime as in Theorem~\uref{genl-results}, fix a subgroup $H \subset \GL_2(\wh{\bZ})$ such that $\Gamma_1(n) \subset H \subset \Gamma_0(n)$, and let $\pi\colon J_H \surjects A$ be the optimal newform quotient that is isogenous to $A_f$. Then
\[
\tst p \nmid \#\p{\f{S_2(\Gamma_0(n),\, \bZ)[e_{f^\perp}]}{\pi^*(H^0(\cA,\, \Omega^1))}}
\]
and $\pi \colon (\cJ_H)_{\bZ_p} \ra \cA_{\bZ_p}$ is smooth, where $\cJ_H$ and $\cA$ denote the N\'{e}ron models over $\bZ$ of $J_H$ and $A$.
%where $S_2(H, \bZ)$ denotes the $\bZ$-module of those weight $2$ cusp forms of level $H$ whose $q$-expansion at ``$\infty$'' lies in $\bZ\llb q \rrb$ and $S_2(H,\, \bZ)[e_{f^\perp}] \ce S_2(H,\, \bZ) \bigcap (S_2(\Gamma_0(n),\, \bQ)[e_{f^\perp}])$ inside $S_2(H,\, \bZ) \tensor_\bZ \bQ$.
\ecor

\bpf
As in \Cref{down}, the multiplicity one theorem supplies a unique isogeny $a$ such that
\be \lab{com-com} \ba 
\xymatrix{
J_{H} \ar@{->>}[r]^-{\pi} \ar@{->>}[d] & A \ar@{->>}[d]^-{a} \\
J_0(n) \ar@{->>}[r]^-{\pi_f} & A_f
}
\ea \ee
commutes. The analogous statement also holds with $J_0(n)$ replaced by $J_{H'}$ for some $H \subset H' \subset \Gamma_0(n)$, so the inclusion 
\[
 \pi^*(H^0(\cA, \Omega^1)) \subset S_2(\Gamma_0(n), \bZ)[e_{f^\perp}]
 \]
 supplied by \cite{CES03}*{Lem.~6.1.6} in the case $H = \Gamma_1(n)$ implies this inclusion for general $H$. Moreover, \eqref{com-com} gives the equality
\[
\tst \# \p{\f{S_2(\Gamma_0(n),\, \bZ)[e_{f^\perp}]}{\pi^*_f(H^0(\cA_f,\, \Omega^1))}} = \#\p{\f{S_2(\Gamma_0(n),\, \bZ)[e_{f^\perp}] }{\pi^*(H^0(\cA,\, \Omega^1))}} \cdot \#\p{ \f{\Lie \cA_f}{ (\Lie a)(\Lie \cA)}},
\]
so the claim follows from \Cref{genl-results}.
\epf

\brem
The proof of \Cref{CES-positive} also shows that the induced map $\cA_{\bZ_p} \xra{a} (\cA_f)_{\bZ_p}$ is \'{e}tale.
\erem

%%%%%%%%%%%%%%%%%%%%%%%%%%%%%%%%%%%%%%

\end{document}